\documentclass[final]{siamltex}
\usepackage{graphicx}
\usepackage{amsmath,amssymb,amsfonts,bm}    
\usepackage{graphics}                           
\usepackage{color}                              
\usepackage{dsfont}
\usepackage{hyperref}                           
\usepackage{multirow}                           
\usepackage{graphicx}                           

\usepackage{makeidx}
\usepackage{subfigure}                          
\usepackage{tikz}
\usepackage{srcltx}                             %
\usepackage[normalem]{ulem} 
\usepackage{xargs,cleveref,enumerate}
\usepackage{srcltx}                             %
\usepackage{fancyhdr}
\usepackage{float}
\usepackage{pgf}
\usetikzlibrary{decorations.pathmorphing}
\usetikzlibrary{fadings}
\usepackage{pgfkeys} 
\usepackage{ifthen}
\usepackage{xargs,cleveref,enumerate}

\usepackage{amsmath,amsxtra,amssymb,latexsym,amsfonts,amscd}
\usepackage{multicol,color}
\usepackage{float}
\usepackage[normalem]{ulem}
\usepackage{pgf,tikz}
\usepackage{mathrsfs}
\usetikzlibrary{arrows}
\usepackage{amssymb}
\usepackage{theorem}
\usepackage{fancyhdr}
\usepackage{enumerate}
\usepackage[margin=1.1in]{geometry}
\usepackage{bbm}
\usepackage[notref,notcite]{showkeys}

\def\disp{\displaystyle}
\def\lm{\lambda}
\def\blm{\bm{\lambda}}
\def\gph{\hbox{}}
\def\gph{\mbox{\rm gph}\,}

\def\O{\Omega}

\def\ph{\varphi}

\def\R{\mathbb{R}}
\def\N{\mathbb{N}}

\def\P{\mathbb{P}}

\def\K{\mathcal{K}}

\def\dom{\mbox{\rm dom}\,}
\def\({\left(}
\def\){\right)}
\def\[{\left[}
\def\]{\right]}
\def\n{\left \|}
\def\en{\right \|}
\def\la{\langle}
\def\ra{\rangle}

\def\ph{\varphi}

\def\gg{\gamma}

\def\ve{\varepsilon}
\def\ref{\mathrm{ref}}
\def\emp{\emptyset}
\def\tto{\rightrightarrows}
\def\gph{\mbox{\rm gph}\,}
\def\dom{\mbox{\rm dom}\,}
\def\epi{\mbox{\rm epi}\,}
\def\dd{\delta}

\def\h{\hfill\triangle}

\def\ox{\bar{x}}
\def\oy{\bar{y}}
\def\oz{\bar{z}}
\def\ov{\bar{v}}

\def\obfy{\bar{\mathbf{y}}}

\def\oR{\Bar{\R}}
\def\al{\alpha}

\def\bfy{\mathbf{y}}

\def\bfY{\mathbf{Y}}
\def\ref{\mathrm{ref}}
\def\obfy{\bar{\mathbf{y}}}
\def\obfY{\bar{\mathbf{Y}}}

\def\bff{\mathbf{f}}

\def\dif{\mathrm{d}}
\def\bfD{\mathbf{D}}

\def\tto{\rightrightarrows}
\def\Lto{\Longrightarrow}
\newcommand{\QVIorig}{\mathbf{QVI}}
\newcommand{\QVIsemi}{\mathrm{QVI}}
\newtheorem{thm}{Theorem}

\newtheorem{cor*}{Corollary}

\newtheorem{defn}{Definition}

\newcommandx{\yaHelper}[2][1=\empty]{%
\ifthenelse{\equal{#1}{\empty}}%
{\ensuremath{ \scriptstyle{ #2 } } } 
{ \raisebox{ #1 }[0pt][0pt]{ \ensuremath{ \scriptstyle{ #2 } } } }  
}\newcommandx{\yrightarrow}[4][1=\empty, 2=\empty, 4=\empty, usedefault=@]{%
\ifthenelse{\equal{#2}{\empty}}
{\xrightarrow{ \protect{ \yaHelper[ #4 ]{ #3 } } } } 
{\xrightarrow[ \protect{ \yaHelper[ #2 ]{ #1 } } ]{ \protect{ \yaHelper[ #4 ]{ #3 } } } }}
\newcommand{\Mosco}{\hspace{-.05cm}\yrightarrow{\scriptscriptstyle \mathrm{M}}[-1pt]\hspace{-.05cm}}

\usepackage[notref,notcite]{showkeys}
\usepackage[textsize=small]{todonotes}
\title{Optimal Control of a Quasi-variational Sweeping Process \thanks{HA, RA, and CNR are partially supported by NSF-DMS 2110263, 2012391,1913004, and Air Force Office of Scientific Research under Award NO: FA9550-19-1-0036 and FA9550-22-1-0248.}}

\author{Harbir Antil\thanks{Department of Mathematical Sciences and the Center for Mathematics and Artificial Intelligence (CMAI), George Mason University, Fairfax, VA 22030, USA ({\tt hantil@gmu.edu},  {\tt tarndt@gmu.edu}, {\tt crautenb@gmu.edu})} \and Rafael Arndt$^\dagger$ \and Boris S. Mordukhovich\thanks{Department of Mathematics, Wayne State University, Detroit, MI 48202 ({\tt aa1086@wayne.edu}). Research of this author was partially supported by the US National Science Foundation under grants DMS-1808978 and DMS-2204519, and by the Australian Research Council under Discovery Project DP-190100555.} \and Dao Nguyen \thanks{Department of Mathematics, University of Michigan, Ann Arbor, MI 48109 ({\tt nntdao@umich.edu}).} \and Carlos N. Rautenberg$^\dagger $}

\begin{document}

\maketitle
\begin{abstract} The paper addresses the study of a class of evolutionary quasi-variational inequalities of the parabolic type arising in the formation and growth models of granular and cohensionless materials. Such models and their mathematical descriptions are highly challenging and require powerful tools of their analysis and implementation. We formulate a space-time continuous optimal control problem for a basic model of this type, develop several regularization and approximation procedures, and establish the existence of optimal solutions {for the time-continuous and space-discrete problem}. Viewing a version of this problem as a controlled quasi-variational sweeping  process leads us to deriving necessary optimality conditions for {the fully discrete problem} by using the advanced machinery of variational analysis and generalized differentiation. 
\end{abstract}\vspace*{-0.1in}

\vspace*{-0.1in}


\section{Introduction}

This paper concerns the optimal selection of a supporting surface for the minimal accumulation of  some granular cohensionless material that is being poured into a known region. The corresponding mathematical model can be formulated as an {\em optimal control problem} for an {\em evolutionary quasi-variational inequality} (QVI), or a {\em quasi-variational sweeping process}, with a gradient type constraint discussed in what follows. The problem is not standard in nature as the control variable acts on the nonconvex constraint set, and thus face significant complexity in establishing {\em well-posedness} and deriving {\em necessary optimality conditions}; see below for more details.

\begin{figure}\label{fig:image}
\includegraphics[scale=0.21]{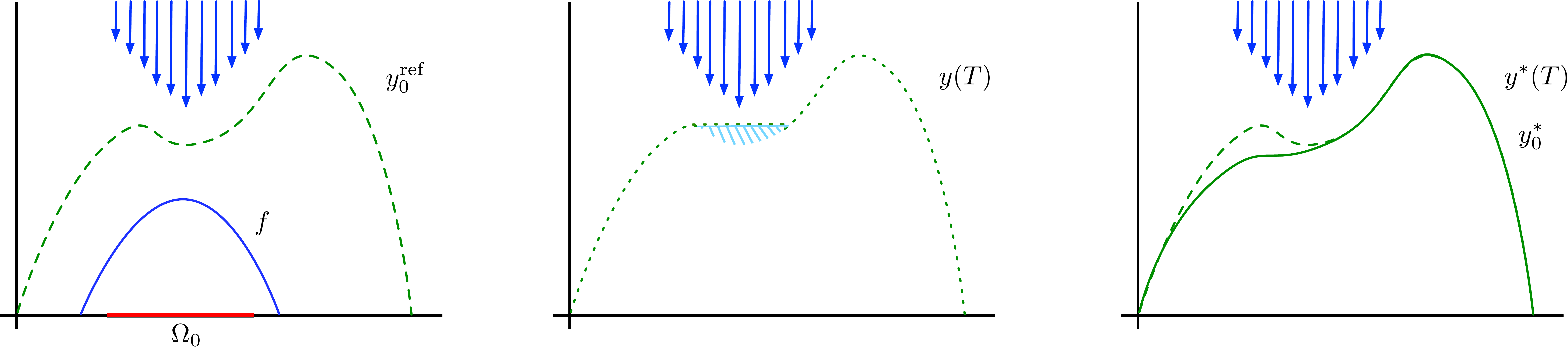}
\caption{(\textbf{LEFT}) Depiction of the initial supporting structure $y_0^{\rm ref}$, the density rate of poured material $f$ and the location of the subdomain $\Omega_0$, where the material should not accumulate. (\textbf{CENTER}) Final resulting shape at time $t=T$ for the material with a very flat angle of repose. (\textbf{RIGHT})  Optimal supporting structure $y_0^*$, which coincides with the final growth shape $y^*$ at time $t=T$ given that no material is accumulating anywhere.}
 \end{figure}

In mathematical terms, the initial supporting surface $y_0$ is defined as a function on a certain domain $\Omega$ {which is} vanishing at the boundary $\partial\Omega$. Suppose that the density rate of the cohensionless granular material that is poured over $y_0^{\mathrm{ref}}$ is known and is denoted by $f$. The resulting final shape of the growth surface is {denoted by} $y$. Furthermore, a subdomain $\Omega_0\subset\Omega$ is provided, where we are supposed to avoid the accumulation of material on a certain time interval $[0,T]$. Assume also that certain perturbations of $y_0^{\mathrm{ref}}$ are allowed while leading us to a shape $y^*_0$, where we aim to maintain the constraints  $0\le y^*-y_0^*\ll y-y_0^{ref}$ over $\Omega_0$ in a prescribed sense. Here $y^*$ is the state corresponding to the initial surface $y_0^*$. A schematic of this behavior in typical cases has been depicted in Figure~\ref{fig:image}.

Having in mind that the problem possesses {insurmountable} difficulties in the original setting (in particular, it can be viewed as the control of the {\em fixed point} of a {\em discontinuous} mapping), we initially tackle a {\em semi-discrete} (in space) version of the problem with a regularized upper bound of the gradient constraint. In this setting, we are able {to prove} {\em existence of feasible solutions} to the resulting {\em QVI} by developing monotone regularization techniques, and then the {\em existence of minimizers} to the overall {\em optimization problem} by properly identifying conditions for the {\em Mosco set convergence} associated to the gradient constraint. {Furthermore}, we develop several regularization and approximation procedures, which allow us to model an {appropriate version} of the basic problem as optimal control of the {\em quasi-variational sweeping process}, which has been never considered in the literature. Nevertheless, applying advanced tools of {\em variational analysis} and {\em second-order generalized differentiation} {enables us to derive} efficient necessary optimality conditions for {fully discretized} quasi-variational sweeping process expressed entirely via given data of the original problem.\vspace*{0.02in}

The rest of the paper is organized as follows. In Section~\ref{sec:ProblemFormulation}, we formulate the {\em original QVI control} problem in appropriate functional spaces {and discuss its regularization.}
The {\em semi-discrete} (in space) QVI is analyzed in Section~\ref{sec:SemiDiscrete}, where the existence and time-regularity are justified. In the same section, the perturbation of solutions with respect to supporting structures is studied. The latter allows us to obtain an existence result for the semi-discrete optimization problem. A formal derivation of stationarity conditions for a regularized problem is provided in Section~\ref{s:formal_deriv}. In Section~\ref{sec:DiscApprox}, we {consider a {\em fully discrete} problem and establish existence of solutions to the corresponding optimization problem.} Section~\ref{sec:coderivative} reviews tools of first-order and second-order variational analysis and generalized differentiation, which allow us to derive {\em necessary optimality conditions} {for the} discretized sweeping control problem with smoothed gradient constraints. The concluding Section~\ref{conclusion} summarizes the major results obtained in this paper with {a discussion on} subsequent {\em numerical} implementations and {a future outlook}.

\section{Problem Formulation and Smoothing}\label{sec:ProblemFormulation} 

Let $y_0:\Omega\to\mathbb{R}$ with $\Omega\subset\mathbb{R}^\mathrm{d}$ be a supporting structure such that $y_0|_{\partial\Omega}=0$, and let $y:(0,T)\times\Omega\to \mathbb{R}$ as $t\in (0,T)$ be the height of the pile of a granular cohensionless material that begins pouring into the domain. Suppose that $y(t)|_{\partial\Omega}=0$, which implies that the material is allowed to abandon the domain freely. The material is characterized by its angle of repose $\theta>0$ that corresponds to the steepest stable angle at which a slope may arise from a point source of the material. The (density)  rate of the granular material being deposited at each point of the domain $\Omega$ is given by $f:(0,T)\times\Omega\to\mathbb{R}$. The mathematical description of such problems was pioneered by Prigozhin and his collaborators \cite{BarrettPrigozhinSandpile,MR3231973, MR3335194,Prigozhin1986,Prigozhin1994,Prigozhin1996,Prigozhin1996a} in the case of homogeneous materials (see also \cite{hr,hr2,hrs}). In this setting, we arrive at the following {\em QVI problem} with respect to the variable $y${, with $p \in [2,\infty]$ and $(\cdot,\cdot)$ denoting the $L^2(\Omega)$ scalar product},\vspace*{0.1in}

\textbf{Problem }($\QVIorig(y_0)$).
Find $y\in L^2(0,T; H_0^1(\Omega))$ with $\partial_ty\in L^2(0,T;H^{-1}(\Omega))$ and
\begin{equation}\label{qvi-state}
y\in\mathcal{K}^p(y,y_0):=\big\{z\in H_0^1(\Omega)\;\big|\;\mbox{ with }\;|\nabla z|_p\le M_p(y,y_0)\big\}
\end{equation}
a.e.\ in $(0,T)$, and for which we have
\begin{equation}\label{qvi-parab}
(\partial_t y-f,v-y)\ge 0
\end{equation}
whenever $v\in\mathcal{K}^p(y,y_0)$ a.e.\ in $(0,T)$. The operator $M_p(w,y_0):\Omega\to\mathbb{R}$ in $\QVIorig(y_0)$ is given by
\begin{equation}\label{qviM}
M_p(w,y_0):=
\begin{cases}
\alpha&\text{if $w>y_0$,} \\
\max\big(\alpha,|\nabla y_0|_{p}\big)&\text{if $w=y_0$,}
 \end{cases}
\end{equation}
where $\alpha:=\tan(\theta)$. In particular, this means that if the material has accumulated, then the gradient constraint is the material dependent one, but if it has not, we may get higher gradients on the supporting surface. This actually permits the material to slide off high slopes into other regions.

A few words are to be said {about the} problem $(\QVIorig(y_0))$. Namely, this {is a} highly {\em nonsmooth} and {\em nonconvex} problem, where the {\em intrinsic  nonconvexity} is induced by the constraint $y\in\mathcal{K}^p(y,y_0)$. Even for the case when $M_p(y,y_0)\equiv C$, a constant, the problem is nonsmooth, and the fact that the gradient is constrained pointwise increases the nonlinearity (with respect to the obstacle constraints) of the overall problem. However, the major difficulty associated to the aforementioned problem is the fact that $M_p$ is {\em discontinuous}, and thus the solution to $\QVIorig(y_0)$ is a {\em fixed point} of a discontinuous mapping. Existence, stability, and the overall analysis for this kind of problems are extremely challenging {and are still open in the most general setting}.

The choice of $p$ in $\QVIorig(y_0)$ determines possible shapes of $y$ and hence the possible structures of the piles. In particular (and formally), if we consider a point source $f$ in the case $p=2$, the structure of $y$ (for a flat $y_0$) corresponds to a growing cone. Other cases like $p=\infty$ would imply that a point source of sand would generate a pyramid structure, where sides are aligned with the horizontal and vertical axis instead. In the latter case, note that $v\in \mathcal{K}^\infty(y,y_0)$ implies that
\begin{equation*}
-M_{\infty}(y,y_0)\le\partial_{x_i}v\le M_{\infty}(y,y_0) \quad \text{ a.e. in } \Omega \text{ and for }i=1,2,\ldots,\mathrm{d}. 
\end{equation*}

A natural {\em optimal control problem} for $\QVIorig(y_0)$ can be described in words as follows. We want to modify $y_0$ (slightly), with respect to some reference structure $y_0^{\mathrm{ref}}$, in order to maintain a certain region of $\Omega_0\subset\Omega$ relatively free of material. This leads us to the following optimization problem with a quasi-variational inequality constraint.\vspace*{0.1in} 

\textbf{Problem} $(\mathbb{P})$. Given $\sigma>0$, $f\in L^2(\Omega)^+$, and $y_0^{\mathrm{ref}}\in L^2(\Omega)$, consider the {\em optimization problem}:
\begin{align*}
&\mathrm{minimize}\qquad \int_0^T\int_{\Omega_0}(y-y_0)\,\mathrm{d}x\,\mathrm{d}t+\frac{\sigma}{2}\int_{\Omega}(y_0-y_0^{\mathrm{ref}})^2 \mathrm{d}x\quad \text{over} \quad y_0 \in H_0^1(\Omega)\\
&\mathrm{subject \:\: to \:\: (s.t.)}\quad  y \text{ solves } \QVIorig(y_0),\\
&\hphantom{\mathrm{subject \:\:to \:\: (s.t.)}}\quad y_0\in \mathcal{A},
\end{align*}
where the constraint set ${\cal A}$ is described by
\begin{equation}\label{qvi-control}
\mathcal{A}:=\big\{z\in H_0^1(\Omega)\;\big|\;y_0^{\mathrm{ref}}+\lambda_0\le z\le y_0^{\mathrm{ref}}+\lambda_1\:\text{ a.e.}\}
\end{equation}
with $\lambda_i\in H_0^1(\Omega)$ for {$i=0,1$} and $\lambda_0\le\lambda_1$ a.e.\ in $\Omega$.
\vspace{.5cm}

{Here $\sigma > 0$ is the given regularization parameter.} Let us discuss some underlying features of the optimization problem $(\mathbb{P})$. This problem can be viewed as an optimal control problem for quasi-variational inequalities with state functions $y\in L^2(0,T;H_0^1(\Omega))$ and control functions $y_0\in H_0^1(\Omega)$ over the {\em parabolic} QVI  \eqref{qvi-parab} subject to the {\em hard/pointwise control constraint} \eqref{qvi-control} and the {\em mixed state-control constraint} \eqref{qvi-state}. This type of optimal control problems are among the most challenging in control theory. As mentioned above, the state-control constraint \eqref{qvi-state} is really complicated from the viewpoint of quasi-variational inequalities. This constraint also creates {trouble to}  handle it from the viewpoint of optimal control {and to derive the} necessary optimality conditions. 

Observe also that we consider $H_0^1(\Omega)$ control perturbations in \eqref{qvi-control} with certain pointwise bounds. In particular, this makes it possible for specific regions to get modified, while other regions of $y_0^{\mathrm{ref}}$ may remain the same. Note {further} that even the application of the direct method of calculus of variations falls short to tackle the existence of optimal {solutions to}  $(\mathbb{P})$. In particular, without certain additional hypotheses, a minimizing sequence $\{y_0^n\}$ of this problem does not allow us to pass from ``$y^n \text{ solves } \QVIorig(y^n_0)$'' to the existence of a function $y^* \text{ solving } \QVIorig(y^*_0)$, where $y^*_0$ is some accumulation point of $\{y_0^n\}$. We discuss the main assumptions needed for the application of the direct method in the next section.\vspace*{0.05in}

{In order to overcome part of these challenges}, we consider certain {\em smoothing approximation procedures} to deal with the mapping $M_p$, which is discontinuous, a major obstacle from both theoretical and numerical viewpoints. To this end, we observe that the mapping $M_p$ can be redefined as
\begin{equation*}
M_ p(w,y_0):=
\begin{cases}
\alpha&\text{if $w>y_0+\epsilon$,}\\
\max\big(\alpha,|\nabla y_0|_{2}\big)\disp\frac{(y_0+\epsilon-w)}{\epsilon}+\disp\alpha\frac{(w-y_0)}{\epsilon}&\text{if $y_0+\epsilon\ge w>y_0$,}\\
\max\big(\alpha,|\nabla y_0|_{2}\big)&\text{if $w=y_0$}.
 \end{cases}
\end{equation*}
Smoother approximations $\Tilde{M}_p$ of $M_p$ may be obtained by using {\em higher-order interpolants} as well as {\em regularizations} of the $\max$ function and the $\mathbb{R}^{\mathrm{d}}$ norm. Due to this, we assume throughout the paper that 
\begin{equation*}
\Tilde{M}_p \text{ is } k \text{ times continuously differentiable for some} \;k\ge 2.
\end{equation*}
The above redefinition of $M_p$ and its approximations allow us to correct a major difficulty associated to the model. Indeed, this induces that the solution $y$ to problem ($\QVIorig(y_0)$) can be equivalently formulated as a {\em fixed point} of a now {\em continuous} mapping, which allows us to employ some perturbation methods.

\section{The semi-discrete problem}\label{sec:SemiDiscrete}
 
In this section, we construct a semi-discrete version of the original QVI control problem $({\cal P})$ involving a {\em space discretization} of the $\Omega$ in the QVI model. Given $\mathbf{f}:(0,T)\to \mathbb{R}^N$, the {\em semi-discrete QVI problem} is formulated as follows: Find $\mathbf{y}:(0,T)\to\mathbb{R}^N$ such that 
\begin{equation}\label{QVI0}\tag{$\QVIsemi_N(\mathbf{y}_0)$}
\mathbf{y}(t)\in \mathcal{K}^p(\mathbf{y}(t),\mathbf{y}_0)\;\mbox{ with }\;\left( \mathbf{y}'(t)-\mathbf{f}(t),\mathbf{v}-\mathbf{y}(t)\right)_{\mathbb{R}^N}\ge 0\;\mbox{ for all }\;\mathbf{v}\in \mathcal{K}^p(\mathbf{y}(t),\mathbf{y}_0)
\end{equation}
and a.e. $t\in(0,T)$. The two most important choices for $\mathcal{K}^p(\mathbf{y}(t),\mathbf{y}_0)$ are $p=2$ and $p=\infty$, where for arbitrary $\mathbf{w}$ and $\mathbf{z}$ 
\begin{align}\label{p2}
&\mathcal{K}^2(\mathbf{w},\mathbf{z}):=\big\{\mathbf{v}\in\mathbb{R}^N\;\big|\;\sqrt{|(\mathbf{D}_1\mathbf{v})_i|^2+|(\mathbf{D}_2\mathbf{v})_i|^2}\le \big(M_2(\mathbf{w},\mathbf{z})\big)_i,\;i=1,\ldots,N\big\},
\end{align}
and
\begin{align}\label{setK}	&\mathcal{K}^\infty(\mathbf{w},\mathbf{z}):=\big\{\mathbf{v}\in\mathbb{R}^N\;\big|\;-
\big(M_\infty(\mathbf{w},\mathbf{z})\big)_i\le(\mathbf{D}_j\mathbf{v})_i\le \big(M_\infty(\mathbf{w},\mathbf{z})\big)_i,\;j=1,\;i=1,2,\ldots,N\big\},
\end{align}
with $\mathbf{D}_1,\mathbf{D}_2\in\mathbb{R}^{N\times N}$ and $M_2,M_{\infty}:\mathbb{R}^N\times\mathbb{R}^N\to\mathbb{R}^N$. Observe that 
$\mathbf{D}_1,$ and $\mathbf{D}_2$ represent discrete approximations of the partial derivatives $\partial/\partial x$ and $\partial/\partial y$, respectively. In this vein, we have that $\mathbf{D}:=(\mathbf{D}_1, \mathbf{D}_2):\mathbb{R}^{N}\to\mathbb{R}^{2N}$ provides an approximation of the gradient.\vspace*{0.05in}

For the rest of possible $p$ values, i.e., $2< p <\infty$,  for arbitrary  $\mathbf{w},\mathbf{z}$ the $\mathcal{K}^p(\mathbf{w},\mathbf{z})$ sets are defined as
\begin{align}\label{Kp}
&\mathcal{K}^p(\mathbf{w},\mathbf{z}):=\big\{\mathbf{v}\in\mathbb{R}^N\;\big|\; (\mathbf{D}\mathbf{v})_i|_p:=\big(|(\mathbf{D}_1\mathbf{v})_i|^p+|(\mathbf{D}_2\mathbf{v})_i|^p)\big)^{\frac{1}{p}}\le\big(M_p(\mathbf{w},\mathbf{z})\big)_i,\quad i=1,\ldots,N\big\},
\end{align}
where the mapping $M_p$ is defined by
\begin{equation}\label{M_p}
\big(M_p(\mathbf{w},\mathbf{z})\big)_i:=
\begin{cases}
\alpha&\text{if $\mathbf{w}_i>\mathbf{z}_i+\epsilon$,} \\
\max\big(\alpha,|(\mathbf{D} \mathbf{z})_i|_{p}\big)\disp\frac{(\mathbf{z}_i+\epsilon-\mathbf{w}_i)}{\epsilon}{\epsilon}+\disp\alpha\frac{(\mathbf{w}_i-\mathbf{z}_i)}{\epsilon}&\text{if $\mathbf{z}_i+\epsilon\ge\mathbf{w}_i>\mathbf{z}_i$,} \\
\max\big(\alpha,|(\mathbf{D}\mathbf{z})_i|_{p}\big)&\text{if $\mathbf{w}_i=\mathbf{z}_i$}
\end{cases}
\end{equation}
with $\mathbf{D}\mathbf{z}:=(\mathbf{D}_1\mathbf{z},\mathbf{D}_2\mathbf{z})$ and $(\mathbf{D}\mathbf{z})_i:=((\mathbf{D}_1\mathbf{z})_i,(\mathbf{D}_2\mathbf{z})_i)$. Although $M_p$ is only continuous, we can consider {\em smooth approximations} $\Tilde{M}_p$ of $M_p$ as explained in the previous section. 


Now we prove that the quasi-variational inequality \eqref{QVI0} admits at least one solution. Although the proof of the following theorem can be inferred from other sources, we include it for the sake of completeness and due to pieces and parts are used later for other arguments.

\begin{thm}{\bf(existence of solutions to semi-discrete QVIs).}\label{thm:existenceQVI} Let $\mathbf{y}_0\in\mathbb{R}^N$, and $\mathbf{f}:(0,T)\to\mathbb{R}^N$ be such that $\mathbf{f}\in L^2(0,T)$. Then there exists a solution $\mathbf{y}:(0,T)\to\mathbb{R}^N$ to \eqref{QVI0} with the properties
\begin{equation}\label{prop}
\mathbf{y}\in C([0,T]) \quad \text{and} \quad \mathbf{y}'\in L^2(0,T). 
\end{equation}
\end{thm}

{\bf Proof}. We split the proof of the theorem into the following four major steps, where each of the steps is of its independent interest.\\[1ex]
{\bf Step~1:} {\em Existence of solutions to the regularized variational inequality.} We confine ourselves to the case where $p=2$, while observing the the general case with $1\le p\le\infty$ can be done similarly. Given $\gamma>0$, consider the nonlinear ordinary differential equation
\begin{equation}\label{ODE1}
\mathbf{y}'(t)=\mathbf{f}(t)-\gamma G\big(t,\mathbf{y}(t)\big),\quad	\mathbf{y}(0)=\mathbf{y}_0
\end{equation}
with the mapping $G$ in the right-hand side of the equation defined by
\begin{equation*}
G(t,\mathbf{y}(t)):=\mathbf{D}^\mathrm{T}(|\mathbf{D}\mathbf{y}(t)|^2_2-M(t)^2)^+\mathbf{D}\mathbf{y}(t),
\end{equation*}
where $M(t):=M_{2}(\mathbf{z}(t),\mathbf{y}_0)$ for an arbitrary $\mathbf{z}\in C(\mathbb{R})$, $\mathbf{y}_0\in\mathbb{R}^N$, and
\begin{equation*}
\big(\mathbf{h}(t)\big)^+:=\big(\max(h_1(t),0),\max(h_2(t),0),\ldots, \max(\tau_M(t),0)\big)
\end{equation*}
for $\mathbf{h}:(0,T)\to\mathbb{R}^N$. Note that the mapping $\mathbb{R}^N\ni\mathbf{h}\mapsto G(t,\mathbf{h})\in\mathbb{R}^N$ is {\em monotone} for each $t$ in $\mathbb{R}^N$, i.e.,
\begin{equation}\label{MonG}
\big\la G(t,\mathbf{h}_1)-G(t,\mathbf{h}_2),\mathbf{h}_1-\mathbf{h}_2\big\ra\ge 0\;\mbox{ whenever }\;\mathbf{h}_1,\mathbf{h}_2\in\mathbf{R}^N,
\end{equation}
and that we have $J(\mathbf{h})'\mathbf{d}=G(t,\mathbf{h})\mathbf{d}$ for the convex function $J(\mathbf{h}):=((|\mathbf{D}\mathbf{h}|^2_2-M(t)^2)^+\mathbf{D}\mathbf{h},\mathbf{D}\mathbf{h})$.

The integral formulation of \eqref{ODE1} is then given by
\begin{equation}\label{ODE2}
\mathbf{y}(t)=\mathbf{y}_0+\int_0^t\mathbf{f}(s)\dif s-\gamma\int_0^tG\big(s,\mathbf{y}(s)\big)\dif s=:\Lambda(\mathbf{y})(s),
\end{equation}
where $\Lambda:C([0,T])\to C([0,T])$ and $\Lambda(\mathbf{y})'(s)\in L^2(0,T)$ for $\mathbf{y}\in C([0,T])$.

To verify the existence of a solution to \eqref{ODE2} for each $\gamma>0$, we use the classical Leray-Schauder theorem. First note that the operator $\Lambda:C([0,T])\to C([0,T])$ is \emph{continuous}. Taking now a sequence $\{\mathbf{y}_n\}$ bounded in $C([0,T])$ ensures that $\{\Lambda(\mathbf{y}_n)\}$ is also bounded in $C([0,T])$, and furthermore $\{\Lambda(\mathbf{y}_n)'\}$ is bounded in $L^2(0,T)$. Indeed, if $C>0$ is such that 
\begin{equation*}
\sup_{n}\|\mathbf{y}_n(t)\|_{C([0,T])}\le C,
\end{equation*}
then we clearly get the estimate
\begin{equation*}
\|\Lambda(\mathbf{y}_n)'\|_{L^2(0,T)}\le\|\mathbf{f}\|_{L^2(0,T)}+\gamma\|\mathbf{D}^\mathrm{T}\|\cdot\|\mathbf{D}\|C T^{1/2}\big(\|\mathbf{D}\|^2 C^2+\|M\|_{C([0,T])}\big).
\end{equation*}
It follows, by the compact embedding of $V:=\{\mathbf{v}\in L^2(0,T)\;|\;\mathbf{v}'\in L^2(0,T)\}$ into $C([0,T])$, that $\Lambda(\mathbf{y}_n)\to\mathbf{g}$ for some $\mathbf{g}\in C([0,T])$ along a subsequence. This tells us that $\Lambda:C([0,T])\to C([0,T])$ is \emph{completely continuous}. Finally in this step, we prove that the set
\begin{equation*}
Y:=\big\{\mathbf{y}\in C([0,T])\;\big|\;\mathbf{y}=\lambda \Lambda(\mathbf{y})\;\text{ for some }\;\lambda\in(0,1)\big\}
\end{equation*}  
is bounded. To this end, observe first that if $\mathbf{y}\in Y$, then $\mathbf{y}(0)=\mathbf{y}_0$ and 
\begin{equation*}
\mathbf{y}'(t)=\lambda\mathbf{f}(t)-\gamma\lambda G\big(t,\mathbf{y}(t)\big).
\end{equation*}
Taking the inner product of $\mathbf{y}$ with the integral from $0$ to $s<T$ gives us
\begin{align*}
\|\mathbf{y}(s)\|^2_2-\|\mathbf{y}_0\|_2^2&=\lambda\int_0^s\mathbf{f}(t)\cdot \mathbf{y}(t)\dif t-\lambda\gamma\int_0^tG(s,\mathbf{y}(s))\mathbf{y}(s)\dif s\\
&\le\lambda\left(\sup_{t\in[0,T]}\|\mathbf{y}(s)\|_2\right)\int_0^T\|\mathbf{f}(s)\|_2 \dif s,
\end{align*}
where we use that $G(s,\mathbf{h})\mathbf{h}\ge 0$ for all $\mathbf{h}\in\mathbb{R}^n$ and all $s\in(0,T)$. It follows that
\begin{equation}
\sup_{t\in [0,T]}\|\mathbf{y}(t)\|_2\le C_1(\mathbf{y}_0,\mathbf{f})<\infty,
\end{equation}
i.e., all elements of $Y$ are bounded. Therefore, the Leray-Schauder theorem yields the existence of a solution $\mathbf{y}^\gamma$ to \eqref{ODE2} for each $\gamma>0$.\\[1ex]
{\bf Step~2:} {\em Uniqueness of solutions to the regularized variational inequality}.
To verify the uniqueness, suppose that we have two solutions $\mathbf{y}_i^\gamma$ for $i=1,2$. Then, since both functions satisfy \eqref{ODE1}, we subtract term by term and test the equation with $\mathbf{y}_1^\gamma-\mathbf{y}_2^\gamma$ with integrating it from $0$ to $s<T$. Thus it follows from the monotonicity in \eqref{MonG} that
\begin{align*}
\|(\mathbf{y}_1^\gamma-\mathbf{y}_2^\gamma)(s)\|^2_2&=-\gamma\int_0^t\big\la G\big(s,\mathbf{y}^\gamma_1(s)\big)-G\big(s,\mathbf{y}_2^\gamma(s)\big),\mathbf{y}_1^\gamma(t)-\mathbf{y}_2^\gamma(t)\big\ra\dif t\le 0,
\end{align*}
which therefore justifies the uniqueness of solutions to \eqref{ODE1}.\\[1ex]
{\bf Step~3:} {\em Existence and uniqueness of solutions to the variational inequality problem.} Arguing similarly to Step~2 allows us to verify the uniform boundedness of solutions $\mathbf{y}_\gamma$ to \eqref{ODE1} with respect to $\gamma>0$. Indeed, we get from \eqref{ODE1} by integrating from $0$ to $t$ and using $G(s,\mathbf{h})\mathbf{h}\ge 0$ for all $\mathbf{h}\in\mathbb{R}^n$ and all $s\in (0,T)$ that
\begin{align*}
\|\mathbf{y}^\gamma(t)\|^2_2-\|\mathbf{y}_0\|_2^2 &\le \left(\sup_{t\in(0,T)}\|\mathbf{y}^\gamma(s)\|_2\right)\int_0^T\|\mathbf{f}(s)\|_2 \dif s.
\end{align*}
This readily implies the estimate
\begin{equation}\label{eq:ybound}
\sup_{\gamma>0}\sup_{s\in[0,T]}\|\mathbf{y}^\gamma(s)\|_2\le C_1(\mathbf{y}_0,\mathbf{f})<\infty.
\end{equation}
By testing in \eqref{ODE1} with an arbitrary $\mathbf{v}\in L^2(0,T)$ such that $\mathbf{v}'\in L^2(0,T)$, we get
\begin{align}\label{Gzero}
&\gamma\int_0^TG\big(s,\mathbf{y}^\gamma(s)\big)\mathbf{v}(s)\dif s=\int_0^T\mathbf{f}(s)\mathbf{v}(s) \dif s-\int_0^T(\mathbf{y}^\gamma)'(s)\mathbf{v}(s)\dif s\\\notag
&\qquad\le\left(\int_0^T\|\mathbf{f}(s)\|^2_2\dif s\right)^{1/2}\left(\int_0^T\|\mathbf{v}(s)\|^2_2\dif s\right)^{1/2}+\int_0^T \mathbf{y}^\gamma(s)\mathbf{v}'(s)\dif s+\mathbf{y}^\gamma(T)\mathbf{v}(T)-\mathbf{y}_0\mathbf{v}(0)
\\\notag
&\qquad\le\left(\int_0^T\|\mathbf{f}(s)\|^2_2 \dif s\right)^{1/2}\left(\int_0^T\|\mathbf{v}(s)\|^2_2 \dif s\right)^{1/2}+C_1(\mathbf{y}_0,\mathbf{f})T^{1/2}\left(\int_0^T\|\mathbf{v}'(s)\|^2_2 \dif s\right)^{1/2}.
\end{align}
Since $V:=\{\mathbf{v}\in L^2(0,T)\;|\;\mathbf{v}'\in L^2(0,T)\}$ is continuously and compactly embedded in $C([0,T])$, we have
\begin{align*}
&\gamma\int_0^TG\big(s,\mathbf{y}(s)\big)\mathbf{v}(s)\le C_2(\mathbf{y}_0,\mathbf{f},T) \left(\left(\int_0^T\|\mathbf{v}(s)\|^2_2\dif s\right)^{1/2}+\left(\int_0^T\|\mathbf{v}'(s)\|^2_2 \dif s\right)^{1/2}\right)
\end{align*}
for some $C_2(\mathbf{y}_0,\mathbf{f},T)$. This yields the estimate
\begin{align*}
\sup_{\gamma>0}\|\gamma G\big(s,\mathbf{y}(s)\big)\|_{V^*}&\le C_2(\mathbf{y}_0,\mathbf{f},T),
\end{align*}
and hence $\sup_{\gamma>0}\|(\mathbf{y}^\gamma)'\|_{V^*}\le C_3(\mathbf{y}_0,\mathbf{f},T)$. 
In particular, we get that 
\begin{align}\label{eq:Derybound}
\sup_{\gamma>0}\|(\mathbf{y}^\gamma)'\|_{L^2(0,T)}&\le C_3(\mathbf{y}_0,\mathbf{f},T).
\end{align}
Note that $\{\mathbf{y}^\gamma\}_{\gamma>0}$ is bounded in $V$, so we can choose a a sequence $\mathbf{y}^n:=\mathbf{y}^{\gamma_n}$ with $\gamma_n\to\infty$ such that $\mathbf{y}^n\rightharpoonup \mathbf{y}^*$ for some $\mathbf{y}^*\in V$. Since $V$ is continuously and compactly embedded in $C([0,T])$, it follows that
\begin{equation*}
\mathbf{y}^n\to \mathbf{y}^* \text{ in } C([0,T])\quad\text{and}\quad (\mathbf{y}^n)'\rightharpoonup (\mathbf{y}^*)' \text{ in } L^2(0,T).
\end{equation*}
Moreover, observe from \eqref{Gzero} that
\begin{align}\label{Gzeros}
&\lim_{n\to\infty}\int_0^TG\big(s,\mathbf{y}^n(s)\big)\mathbf{v}(s)\dif s=\int_0^T\big((|\mathbf{D}\mathbf{y}^*(t)|^2_2-M(t)^2)^+\mathbf{D}\mathbf{y}^*(t),\mathbf{D}\mathbf{v}(s)\big)\dif s=0,
\end{align}
from which we deduce that $|\mathbf{D}\mathbf{y}^*(t)|_2\le M(t)$, i.e., $\mathbf{y}^*\in\mathcal{K}^2(\mathbf{z}(t),
\mathbf{y}_0)$.
Testing further \eqref{ODE1} with $\mathbf{w}=\mathbf{v}-\mathbf{y}^n$ as $\mathbf{v}\in\mathcal{K}^2(\mathbf{z}(t),\mathbf{y}_0)$ gives us the equality
\begin{equation}\label{ODEweak}
\begin{split}
\int_0^T\big\la(\mathbf{y}^n)'(t)-\mathbf{f}(t),\mathbf{v}(t)-\mathbf{y}^n(t)\big\ra\dif t&=\gamma\int_0^T\big\la G\big(t,\mathbf{v}(t)\big)-G\big(t,\mathbf{y}(t)\big),\mathbf{v}(t)-\mathbf{y}(t)^n\big\ra,
\end{split}
\end{equation}
where the condition $G(t,\mathbf{v})=0$ is used. Employing the fact that $\mathbf{h}\mapsto G(t,\mathbf{h})$ is monotone, we have that the right hand-side of \eqref{ODEweak} is nonnegative. Passing there to the limit as $n\to\infty$ leads us to
\begin{equation}
\begin{split}\label{VIweak}	\int_0^T\big\la (\mathbf{y}^*)'(t)-\mathbf{f}(t),\mathbf{v}(t)-\mathbf{y}^*(t)\big\ra\dif t&\ge 0.
\end{split}
\end{equation}
Since $\mathbf{v}$ was chosen arbitrary, a simple density device shows that $\mathbf{y}^*$ solves the variational inequality
\begin{equation}\label{VI0}
\mathbf{y}(t)\in\mathcal{K}^2\big(\mathbf{z}(t),\mathbf{y}_0\big)\;\big|\;\left\la \mathbf{y}'(t)-\mathbf{f}(t),\mathbf{v}-\mathbf{y}(t)\right\ra_{\mathbb{R}^N}\ge 0\;\mbox{ for all }\;\mathbf{v}\in\mathcal{K}^2\big(\mathbf{z}(t),\mathbf{y}_0\big),
\end{equation}
and the claimed uniqueness follows by monotonicity arguments.\\[1ex]
{\bf Step~4:} {\em Existence of solutions to the quasi-variational inequality problem.} Denote by $\mathbf{y}=S(\mathbf{z})$ the (single-valued by Step~3) solution mapping of the variational inequality \eqref{VI0}. Arguing similarly to Step~3 ensures that the mapping $S:V\to V$ is compact. Furthermore, by the estimate
\begin{equation*}
M_2\big(\mathbf{z}(t),\mathbf{y}_0\big)\le\max(\alpha,|\mathbf{D}\mathbf{y}_0|_{\infty})=:\beta
\end{equation*}
we deduce that $S$ maps $\mathcal{K}^2_\beta$ into $\mathcal{K}^2_\beta$, where
\begin{equation*}
\mathcal{K}^2_\beta:=\big\{\mathbf{v}\in C([0,T])\;\big|\;|\mathbf{D}\mathbf{y}(t)|_2\le\beta\;\text { a.e.}\big\}.
\end{equation*}
Employing finally Schauder's fixed point theorem yields the existence of a fixed point $\mathbf{y}=S(\mathbf{y})$, and therefore the quasi-variational inequality \eqref{QVI0} admits a solution satisfying \eqref{prop}. This verifies the statement of Step~4 and thus completes the proof of the theorem. $\h$\vspace*{0.08in}

Now we formulate the following {\em optimal control problem with the \eqref{QVI0} constraints}. The previous theorem allows us to pose the problem in a slightly more regular space than chosen initially.\vspace*{0.1in}

\textbf{Problem} $(\mathbb{P}_N)$. Given a number $\sigma>0$, a nonnegative (i.e., with nonnegative components) mapping $\mathbf{f}:(0,T)\to \mathbb{R}^N$, and vectors $\mathbf{a},\mathbf{y}_0^{\mathrm{ref}}\in \mathbb{R}^N$, consider the following optimal control problem for \eqref{QVI0}:
\begin{align*}
&\mathrm{minimize} \qquad J(\mathbf{y},\mathbf{y}_0):=\int_0^T\big\la\mathbf{a}, \mathbf{y}(t)-\mathbf{y}_0\big\ra\mathrm{d}t+\frac{\sigma}{2}\big\la\mathbf{y}_0-\mathbf{y}_0^{\mathrm{ref}},\mathbf{y}_0-\mathbf{y}_0^{\mathrm{ref}}\big\ra\quad \text{over}\quad \mathbf{y}_0 \in \mathbb{R}^ N\\
&\mathrm{subject \:\: to \:\:}\quad \mathbf{y} \text{ solves } \QVIsemi(\mathbf{y}_0),\\	&\hphantom{\mathrm{subject \:\: to \:\:}\quad}  \mathbf{y}\in V:=\big\{\mathbf{v}\in  L^2(0,T)\;\big|\;\mathbf{v}'\in L^2(0,T)\big\},\\	&\hphantom{\mathrm{subject \:\: to \:\:}\quad}  \mathbf{y}_0\in \mathcal{A},
\end{align*}
where the latter {\em control constraint} set is defined by
\begin{equation*}
\mathcal{A}:=\big\{\mathbf{z}\in \mathbb{R}^ N\;\big|\;\mathbf{y}_0^{\mathrm{ref}}+\bm{\lambda}_0\le \mathbf{z}\le \mathbf{y}_0^{\mathrm{ref}}+\bm{\lambda}_1\big\},
\end{equation*}
with $\bm{\lambda}_0,\bm{\lambda}_1\in \mathbb{R}^ N$ such that $0\le \bm{\lambda}_0\leq\bm{\lambda}_1$.\vspace*{0.12in}

Our next goal is to verify the {\em existence of solutions} to the formulated optimal control problem $(\mathbb{P}_N)$. Before this, recall the notion of {\em Mosco convergence} for sets in reflexive Banach spaces. \vspace*{-0.05in}

\begin{defn}{\bf(Mosco convergence).}\label{definition:MoscoConvergence}
Let $\K$ and $\K_n$ as $n\in\mathbb{N}$ be nonempty, closed, and convex subsets of a reflexive Banach space $V$. Then the sequence $\{\K_n\}$ is said to converge to $\K$ in the sense of Mosco as $n\rightarrow\infty$, which is signified by $$\K_n\Mosco\K,$$ if the following two conditions are satisfied:
\begin{enumerate}[\upshape(I)]
\item\label{itm:1}  For each $w\in \K$, there exists $\{w_{n'}\}$ such that $w_{n'}\in \K_{n'}$ for $n'\in \mathbb{N}'\subset\mathbb{N}$ and $w_{n'}\rightarrow w$ in $V$.
\item\label{itm:2} If $w_n\in \K_n$ and $w_n\rightharpoonup w$ in $V$ along a subsequence, then $w\in \K$.
\end{enumerate}
\end{defn}\vspace*{0.05in}

Here is the aforementioned existence theorem for the formulated optimal control problem. \vspace*{-0.05in}

\begin{thm}{\bf(existence of optimal solutions to \eqref{QVI0}).}\label{thm:ExistPN} The optimal control problem $(\mathbb{P}_N)$ for \eqref{QVI0} admits an optimal solution.
\end{thm}
{\bf Proof}. We spit the proof of the theorem into the two major steps.\\[1ex]
{\bf Step~1}: {\em Properties of minimizing sequences in $(\mathbb{P}_N)$}. Observe first that Theorem~\ref{thm:existenceQVI} tells us that for each $\mathbf{y}_0\in\mathcal{A}$ there exists a $\mathbf{y}\in V$ solving $\QVIsemi_N(\mathbf{y}_0)$. This yields the existence of a minimizing sequence $\{(\mathbf{y}_n,\mathbf{y}_0^n)\}$ for problem $(\mathbb{P}_N)$, i.e., for each $n\in\mathbb{N}$ we have
\begin{equation*}
(\mathbf{y}_n,\mathbf{y}_0^n)\in V\times\mathcal{A}, \quad \mathbf{y}_n \text{ solves }\QVIsemi(\mathbf{y}_0^n) \quad \text{ with } J(\mathbf{y}_n,\mathbf{y}_0^n)\to \inf J\;\mbox{ as }\;n\to\infty.
\end{equation*}
Since $\mathbf{y}_0^n\in\mathcal{A}$ for all $n\in\mathbb{N}$, this implies that for every $n\in\mathbb{N}$ there exists a subsequence of the minimizing sequence (no relabeling) and $\mathbf{y}_0^*\in\mathcal{A}$ such that
\begin{equation*}
\mathbf{y}_0^n\to\mathbf{y}_0^*\;\mbox{ as }\;n\to\infty.
\end{equation*}
Taking into account that the solutions $\mathbf{y}_n$ of $\QVIsemi(\mathbf{y}^n_0)$ are in $V$ and deducing from \eqref{eq:ybound} and \eqref{eq:Derybound} that
\begin{equation}\label{eq:boundsyn}
\sup_{s\in [0,T]}\|\mathbf{y}_n (s)\|_2\le \sup_{\bfy_0\in\mathcal{A}}C_1(\mathbf{y}_0,\mathbf{f})<\infty \qquad \text{ and }\qquad \|\mathbf{y}'_n\|_{L^2(0,T)}\le \sup_{\bfy_0\in\mathcal{A}}C_3(\mathbf{y}_0,\mathbf{f},T)<\infty
\end{equation}
with $C_1(\mathbf{y}_0,\mathbf{f})$ and $C_3(\mathbf{y}_0,\mathbf{f},T)$ being independent of $n$, let us check that these bounds are uniform in $n\in\mathbb{N}$. To verify the uniformity, we get from $\mathbf{y}_n\in C([0,T])$ and the proof of Theorem~\ref{thm:existenceQVI} that for each $n\in\mathbb{N}$ there exists  $\mathbf{z}_k\in V$ satisfying the equation
\begin{equation}\label{ODE3}
\mathbf{z}_k(t)=\mathbf{y}^n_0+\int_0^t\mathbf{f}(s)\dif s-k\int_0^tG^n\big(s,\mathbf{z}_k(s)\big)\dif s
\end{equation}
where the integrand $G^n(s,\mathbf{z}(s))$ is given by
\begin{equation*}
G^n(s,\mathbf{z}(s)):=\mathbf{D}^*(|\mathbf{D}\mathbf{z}(t)|^2_2-M^n(t)^2)^+\mathbf{D}\mathbf{z}(t)
\end{equation*}
with $D^*$ standing for the matrix transposition/adjoint operator, and with the mapping $M^n$ defined by
\begin{equation}\label{eq:Mn}
M^n(t):=M_2\big(\mathbf{y}_n(t),\mathbf{y}_0^n\big).
\end{equation}
Observe further by the proof of Theorem~\ref{thm:existenceQVI} that we have the convergence
\begin{equation*}
\mathbf{z}_k\to\mathbf{y}_n\text{ in } C([0,T]) \qquad \text{and}\qquad (\mathbf{z}_k)'\rightharpoonup (\mathbf{y}_n)' \text{ in } L^2(0,T),
\end{equation*}
and that the following bounds are satisfied:
\begin{equation*}
\sup_{k\in\mathbb{N}}\sup_{s\in[0,T]}\|\mathbf{z}_k (s)\|_2\le C_1(\mathbf{y}_0,\mathbf{f}) \qquad \text{ and }\qquad  \sup_{k\in \mathbb{N}}\|\mathbf{z}'_k\|_{L^2(0,T)}\le C_3(\mathbf{y}_0,\mathbf{f},T),
\end{equation*}
This verifies \eqref{eq:boundsyn} by noting that $\sup_{\bfy\in\mathcal{A}}C_1(\mathbf{y}_0,\mathbf{f})$ and $\sup_{\bfy\in\mathcal{A}}C_3(\mathbf{y}_0,\mathbf{f},T)$ are finite.

It follows from \eqref{eq:boundsyn} that, along a subsequence (no relabeling), we have
\begin{equation}\label{eq:yn_convergence}
\mathbf{y}_n\to \mathbf{y}^* \text{ in } C([0,T]) \qquad \text{and}\qquad (\mathbf{y}_n)'\rightharpoonup (\mathbf{y}^*)' \text{ in } L^2(0,T)
\end{equation}
for some $\mathbf{y}^*\in V$, which is an optimal solution to $(\mathbb{P}_N)$ as shown below.\\[1ex]
{\bf Step~2:} {\em The limiting function $\mathbf{y}^*$ is a solution to the quasi-variational inequality $\QVIsemi(\mathbf{y}^*_0)$}. It follows from \eqref{eq:yn_convergence} that the mapping $M^n$ defined in \eqref{eq:Mn} is such that
\begin{equation*}
M^n\to M^* \text{ in } C([0,T]) \quad\text{ with } 
\quad M^*(t):=M_2\big(\mathbf{y}_0^\ast,\mathbf{y}^*(t)\big)
\end{equation*}
and that $M^n(t)\ge\alpha>0$ for all $n$ by definition. Thus we now show that the convergence
\begin{equation}\label{eq:Mosqui}
\mathscr{K}^2(\mathbf{y}_n,\mathbf{y}_0) \Mosco \mathscr{K}^2(\mathbf{y}^*,\mathbf{y}_0)
\end{equation} 
in the sense of Mosco in the $V$ topology holds true, where
\begin{equation*}
\mathscr{K}^2(\mathbf{z},\mathbf{y}_0):=\big\{\mathbf{w}\in V\;\big|\;\mathbf{w}(t)\in \mathcal{K}^2(\mathbf{z}(t),\mathbf{y}_0) \text{ for all } t\in [0,T]\big\}.
\end{equation*}
This clearly follows for item \eqref{itm:2} in  Definition~\ref{definition:MoscoConvergence}: If $\mathbf{w}_n\in \mathscr{K}^2(\mathbf{y}_n,\mathbf{y}_0^n)$ and $\mathbf{w}_n\rightharpoonup \mathbf{w}^\ast$ in $V$ for some $\mathbf{w}^\ast$, then $\mathbf{w}^\ast\in \mathscr{K}^2(\mathbf{y}^*,\mathbf{y}_0^*)$. Indeed, since $V$ is continuously and compactly embedded in $C([0,T])$, we observe that $\mathbf{w}_n\to \mathbf{w}^\ast$  in  $C([0,T])$. Employing then the estimate
\begin{equation*}
\sqrt{|(\mathbf{D}_1\mathbf{w}_n(t))_i|^2+|(\mathbf{D}_2\mathbf{w}_n(t))_i|^2}\le\big(M_2(\mathbf{y}_n(t),\mathbf{y}_0^n)\big)_i
\end{equation*}
for $t\in[0,t]$ and $i=1,\ldots,N$ tells us that
\begin{equation*}
\sqrt{|(\mathbf{D}_1\mathbf{w}^\ast)_i|^2+|(\mathbf{D}_2\mathbf{w}^\ast)_i|^2}\le \big(M_2(\mathbf{y}^*,\mathbf{y}_0^*)\big)_i,
\end{equation*}
i.e., $\mathbf{w}(t)\in \mathcal{K}^2(\mathbf{y}^*(t),\mathbf{y}^*_0)$ for all $t\in [0,T]$, which thus verifies the statement.

Now we turn the attention to \eqref{itm:1} in Definition \ref{definition:MoscoConvergence}. Note that $M^n\ge\alpha>0$ and $M^n\to M^*$ in $C([0,T])$, and so the positive numbers
\begin{equation*}
\beta_n:=\left(1+\frac{\|M^n-M^*\|_{C([0,T])}}{\alpha}\right)^{-1}
\end{equation*}
are such that $\beta_n\uparrow 1$, and that for $\mathbf{w}^*\in \mathscr{K}^2(\mathbf{y}^*,\mathbf{y}^*_0)$ we have $\beta_n\mathbf{w}^*\in \mathscr{K}^2(\mathbf{y}_n,\mathbf{y}^n_0)$ and $\beta_n\mathbf{w}^*\to \mathbf{w}^*$ in $V$ as $n\to\infty$. This therefore verifies \eqref{eq:Mosqui}. Hence the set convergence in \eqref{eq:Mosqui} implies that the function  $\mathbf{y}^*\in \mathscr{K}^2(\mathbf{y}^*,\mathbf{y}^*_0)$ satisfies the inequality
\begin{equation*}
\begin{split}
\int_0^T\big((\mathbf{y}^*)'(t)-\mathbf{f}(t),\mathbf{v}(t)-\mathbf{y}^*(t)\big)\dif t&\geq 0 \quad \text{ for all }\mathbf{v}\in \mathscr{K}^2(\mathbf{y}^*,\mathbf{y}^\ast_0).
\end{split}
\end{equation*}	
Employing the standard density arguments shows that $\mathbf{y}^*$ is actually a solution to the quasi-variational inequalities $\QVIsemi(\mathbf{y}^*_0)$ while justifying in this way the statement of Step~2.

Finally, the lower semicontinuity of the objective functional ensures that
\begin{equation*}
J(\mathbf{y}^*,\mathbf{y}^*_0)\le\liminf_{n\to\infty} J(\mathbf{y}^n,\mathbf{y}^n_0)=\lim_{n\to\infty} J(\mathbf{y}^n,\mathbf{y}^n_0)=\inf J,
\end{equation*}
which thus completes the proof of the theorem. $\h$\vspace*{0.05in}

{In the above result we have shown existence of solution to $(\mathbb{P}_N)$. Before, we introduce the fully discrete problem and provide a rigorous derivation of the first order optimality conditions, we consider a formal derivation of the first order stationarity conditions for a regularized version of $(\mathbb{P}_N)$. The aim of this upcoming section is give a flavor of the first order conditions and provide a potential alternative to numerically solve $(\mathbb{P}_N)$.}


\section{Regularized Problem and Stationarity Conditions}\label{s:formal_deriv}

The following {\em regularized problem} is obtained from problem $(\mathbb{P}_N)$ by a natural regularization of its {\em quasi-variational constraint} {(see \eqref{ODE1})}:\vspace*{0.05in}

\textbf{Problem} $(\widetilde{\mathbb{P}}_N)$. Given numbers $\sigma,\gamma>0$, a mapping $\mathbf{f}:(0,T)\to \mathbb{R}^N$ with nonnegative components, and vectors $\mathbf{a},\mathbf{y}_0^{\mathrm{ref}}\in\mathbb{R}^ N$, consider the regularized problem
\begin{align*}
&\mathrm{minimize} \qquad J(\mathbf{y},\mathbf{y}_0):=\int_0^T\big\la \mathbf{a},(\mathbf{y}(t)-\mathbf{y}_0)\big\ra\mathrm{d}t
+\frac{\sigma}{2} | \mathbf{y}_0-\mathbf{y}_0^{\mathrm{ref}} |_2^2 \quad \text{over} \quad \mathbf{y}_0\in\mathbb{R}^ N
\end{align*}
subject to $\mathbf{y}\in V$ solving the {\em primal state equation}
\begin{equation}
\begin{split}\label{ODE11}
\mathbf{y}'(t)&=\mathbf{f}(t)-\gamma G(t,\mathbf{y}(t), \mathbf{y}_0),\\
\mathbf{y}(0)&=\mathbf{y}_0
\end{split}
\end{equation}
with $G(t,\mathbf{y}(t),\mathbf{y}_0):={\mathbf{D}^\mathrm{T}}\max_\epsilon\left(0,|\mathbf{D}\mathbf{y}(t)|^2_2-\tilde{M}_p(\mathbf{y}(t),\mathbf{y}_0)^2\right)\mathbf{D}\mathbf{y}(t)$ where $\max_\epsilon$ is a smooth approximation of the $\max$ operator, 
and 
\begin{equation*}
\mathbf{y}_0 \in \mathcal{A}:=\big\{\mathbf{z}\in \mathbb{R}^ N\;\big|\;\mathbf{y}_0^{\mathrm{ref}}+\bm{\lambda}_0\le \mathbf{z}\le\mathbf{y}_0^{\mathrm{ref}}+\bm{\lambda}_1\big\},
\end{equation*}
where $\bm{\lambda}_0,\bm{\lambda}_1\in \mathbb{R}^N$ are such that $0\le\bm{\lambda}_0\le\bm{\lambda}_1$.\vspace*{0.05in}

Let us provide a formal derivation of {\em stationarity conditions} for the above regularized problem by using the  {\em Lagrangian formalism}. To proceed, we introduce the {\em Lagrangian functional} 
\begin{align*}
&\mathcal{L}(\mathbf{y},\mathbf{y}_0,\mathbf{p})
=J(\mathbf{y},\mathbf{y}_0)-\left(\int_0^T\Big\la\mathbf{p}(t),\Big(\mathbf{y}'(t)+\gamma  G(t,\mathbf{y}(t),\mathbf{y}_0)-\mathbf{f}(t)\Big)\Big\ra\dif t\right)
\end{align*}
and observe that a variation of $\mathcal{L}$ with respect to $\mathbf{y}$ at a {\em stationary point} $(\mathbf{y},\mathbf{y}_0,\mathbf{p})$ leads us to the state equation \eqref{ODE11}. Applying further integration by parts to the term $\int_0^T\big\la \mathbf{p}(t),\mathbf{y}'(t)\big\ra\dif t$,
we arrive at
\begin{align*}
\begin{array}{ll}
&\mathcal{L}(\mathbf{y},\mathbf{y}_0,\mathbf{p})=J(\mathbf{y},\mathbf{y}_0)\\
&-\disp\left(\int_0^T \Big(-\big\la\mathbf{y}(t),\mathbf{p}'(t)\big\ra
+\gamma\big\la\mathbf{p}(t),G(t,\mathbf{y}(t),\mathbf{y}_0)\big\ra-\big\la\mathbf{p}(t),\mathbf{f}(t\big\ra)\Big)\dif t +\big\la\mathbf{p}(T),\mathbf{y}(T)\big\ra-\big\la\mathbf{p}(0), \mathbf{y}(0)\big\ra\right).
\end{array}
\end{align*}
To derive the adjoint system, compute a variation of $\mathcal{L}$ with respect to $\mathbf{y}$ at
the stationary point $(\mathbf{y},\mathbf{y}_0,\mathbf{p})$ in the direction $\mathbf{h}$ and get in this way the relationships:
\begin{equation*}
\begin{array}{ll}
0=\mathcal{L}_{\mathbf{y}}(\mathbf{y},\mathbf{y}_0,\mathbf{p})(\mathbf{h})=
\disp\int_0^T\big\la\mathbf{h}(t),\mathbf{a}\big\ra\dif t-\disp\Bigg(&\disp\int_0^T\Big(-\big\la\mathbf{h}(t),\mathbf{p}'(t)\big\ra+\gamma\big\la\mathbf{p}(t),\disp G_{\mathbf{y}}(t,\mathbf{y}(t),\mathbf{y}_0)\mathbf{h}(t)\big\ra\Big)\dif t\\
&+\disp\big\la\mathbf{p}(T),\mathbf{h}(T)\big\ra\Bigg),
\end{array}
\end{equation*}
where we use that $\mathbf{y}_0$ is fixed and thus its variation is equal to zero. Choosing first that $\mathbf{h}$ to be compactly supported and then considering the general case brings us to the following {\em adjoint equation} and its {\em boundary condition}: Find 
$\mathbf{p}$ solving the {adjoint} system
\begin{equation}\label{eq:adjcont}
\begin{cases}-\mathbf{p}'(t)
+\gamma\big\la G_{\mathbf{y}}(t,\mathbf{y}(t),\mathbf{y}_0), \mathbf{p}(t)\big\ra=\mathbf{a},\quad t\in(0,T),\\
\mathbf{p}(T)=0. 
\end{cases}	
\end{equation}
Finally, the minimization of $\mathcal{L}$ with respect to $\mathbf{y}_0$ and subject to $\mathbf{y}_0\in\mathcal{A}$ leads us to the variational inequality for the control variable $\mathbf{y}_0$ formulated as follows:
\begin{equation}\label{eq:VIcont}
\left\langle\sigma(\mathbf{y}_0-\mathbf{y}_0^{\rm ref}), \widehat{\mathbf{y}}-\mathbf{y}_0\right\rangle
-\gamma\int_0^T\big\langle\mathbf{p}(t), G_{\mathbf{y}_0}(t,\mathbf{y}(t),\mathbf{y}_0)
(\widehat{\mathbf{y}}-\mathbf{y}_0)\big\rangle\ge 0  \quad 
\mbox{for all }\;\widehat{\mathbf{y}}\in\mathcal{A}. 
\end{equation}
To summarize, the stationarity system corresponding to the above regularized problem 
is given by the relationships\eqref{ODE11}, \eqref{eq:adjcont}, and \eqref{eq:VIcont}.\vspace*{0.05in}

\section{Quasi-Variational Sweeping Process and Discrete Approximations}\label{sec:DiscApprox}\vspace*{-0.05in}

First we recall the construction of the {\em normal cone} to a convex set $\Theta$ at a point $x$ defined by 
\begin{equation}\label{nor-cone}
N_\Theta(\bar x):=\left\{\begin{array}{ll}
{\big\{ x^* \, : \, }  \la x^*,x-\ox\ra\le 0\;\mbox{ for all }\;x\in\Theta\big\}&\mbox{if }\;{\ox \in\Theta},\\
\emp&\mbox{otherwise}.
\end{array}\right.
\end{equation}
Therefore, the convexity of the sets ${\cal K}^p(y,y_0)$ from \eqref{qvi-state} allows us to rewrite the semi-discrete quasi-variational inequality problem from Section~3 in the form
as a {\em quasi-variational sweeping process}
\begin{equation}\label{eq:DiscQVI}\tag{$\QVIsemi_N(\mathbf{y}_0)$}-\mathbf{y}'(t)\in F\big(\bfy(t),\bfy_0\big):= N_{\mathcal{K}^p(\mathbf{y}(t),\mathbf{y}_0)}\big(\bfy(t)\big)-\mathbf{f}(t).
\end{equation}
Note that the classical (uncontrolled) sweeping process was introduced by Moreau in the 1970s motivated by applications to elastoplasticity; see \cite{mor_frict} with the references to his original publications. A characteristic feature of Moreau's sweeping process and its modifications is that the moving set under the normal cone operator depends on time in a certain continuous way. We refer the reader to the excellent recent survey in \cite{bt} with the comprehensive bibliography therein concerning various theoretic aspects and many applications of Moreau's sweeping process and its further extensions.

Since the Cauchy problem for the aforementioned sweeping processes admits a {\em unique solution} due to the maximal monotonicity of the normal cone operator \cite{bt}, the consideration of any optimization problem for such processes is out of question. This is quite opposite to optimal control theory for {\em Lipschitzian} differential inclusions of the type $\dot x\in F(x)$ and the classical theory for systems governed by differential equations $\dot x=f(x,u),\;u\in U$, and their PDE counterparts.

Starting with \cite{chhm1}, various optimal control models for sweeping dynamics have been formulated rather recently {including derivation of optimality conditions.}
They include: problems with moving sets depending on time and control variables \cite{chhm1,chhm2}, problems with controls in associated ODEs \cite{bk}, problems with controls in additive perturbations of the dynamics \cite{ac,pfs,zeidan}, problems with controls in both moving sets and dynamics \cite{ccmn,cm3}. The cited papers impose different assumptions on the problem data, develop diverse approximation techniques, derive various sets of necessary optimality conditions, and contain references to other publications in these directions. But the common point of all these models for controlled sweeping processes is a {\em highly non-Lipschitzian} (in fact, discontinuous) nature of the sweeping dynamics, which restricts the usage of variational machinery employed in the study of Lipschitzian differential inclusions. Observe also that the very definition of the normal cone \eqref{nor-cone} and their nonconvex extensions yields the unavoidable presence of {\em pointwise state} and {\em mixed state-control} constraints of {\em irregular} types, which are among the most challenging issues even in classical theory.

{Having said that, }
we emphasize that---to the best of our knowledge---no optimal control problems have been considered for sweeping processes with moving sets depending not only on time and control variables but on {\em state variables} as well, which is the essence of {\em quasi-variational} vs.\ variational inequalities. This is the case of the \eqref{eq:DiscQVI}  and \eqref{eq:QVIdd} problems studied in what follows.  

Our approach is based on the {\em method of discrete approximations} and tools of generalized differentiation developed in \cite{m95} to derive necessary optimality conditions in optimal control problems for Lipschitzian differential inclusions with finite-dimensional state spaces and then extended in \cite[Chapter~6]{m-book} to infinite-dimensional systems. Since the Lipschitz continuity is crucial in the device of \cite{m95,m-book} and related publications, the extension of this method to the non-Lipschitzian sweeping dynamics requires significant improvements, 
{this has} been accomplished in \cite{ccmn,cm3,chhm1,chhm2} and other papers for different type of controlled sweeping processes associated with variational inequalities. Here we develop some aspects of this method for optimal control of the quasi-variational sweeping process under consideration.\vspace*{0.03in}   

According to the general scheme of the discrete approximation method, we introduce now the {\em fully discretized} (in time and space) form of the quasi-variational inequality \eqref{QVI0} by using for simplicity the {\em uniform Euler scheme} in the replacement of the time derivative $\dot x$ by finite differences. For this matter, 
take any natural number $M\in\N$ and consider the {\em discrete grid/mesh} on $(0,T)$ defined by
\begin{equation*}
T_M:=\big\{0,\tau_M,\ldots,T-\tau_M,T\big\},\quad \tau_M:=\dfrac{T}{M},
\end{equation*}
with the {\em stepsize of discretization} {$\tau_M$} and the {\em mesh points} $t^M_j:=j\tau_M$ as $j=0,\ldots,M$. Then the quasi-variational inequality in \eqref{QVI0} is replaced by
\begin{equation}\label{eq:QVIdd}\tag{$\mathrm{QVI}_N^M(\mathbf{y}_0)$}
\mathbf{y}_{j}^M\in \mathcal{K}^p(\mathbf{y}_0,\mathbf{y}_{j}^M) \quad\Bigg|\quad\left(\frac{\mathbf{y}_{j}^M-\mathbf{y}_{j-1}^M}{\tau_M}-\mathbf{f}_j^M,\mathbf{v}-\mathbf{y}_j^M\right)_{\mathbb{R}^N}\geq 0\;\mbox{ for all }\;\mathbf{v}\in \mathcal{K}^p(\mathbf{y}_0,\mathbf{y}_{j}^M)
\end{equation}
with the discrete time $j=1,\ldots,M$ and the rate discretization
\begin{equation}\label{f-discr}
\mathbf{f}_j^M=\int_{(j-1)\tau_M}^{j\tau_M}\mathbf{f}(t)\mathrm{d}t \qquad j=1,\ldots,M.
\end{equation}
Equivalently, \eqref{eq:QVIdd} can be written as the {\em discretized quasi-variational sweeping process} 
\begin{equation}\label{e:dis-incl2}
\bfy_j^M\in\bfy_{j-1}^M+\tau_M
F_j^M(\bfy_j^M,\bfy_0),\quad j=1,\ldots,M,
\end{equation}
where the feasible discrete velocity mappings $F_j^M$ are defined by
\begin{equation}\label{mapF1}
F_j^M(\mathbf{y},\bfy_0):=-N_{\K^{p}(\bfy,\bfy_0)}(\bfy)+\bff_j^M,\quad j=1,\ldots,M,
\end{equation}
via the normal cone operator of the state and control dependent set $\K^{p}(\bfy,\bfy_0)$.\vspace*{0.05in}

Given $\{\bfy_j^M\}$ satisfying \eqref{eq:QVIdd}, its {\em piecewise linear extension} $\bfy^M(t)$ to the
continuous-time interval $(0,T)$, i.e., the {\em Euler broken line}, is defined by 
\begin{equation*}
\bfy^M(t):=\sum_{j=1}^M\bfy_j^M\chi_{I_j}(t), \qquad \text{where}\quad I_j=\big[(j-1)\tau_M,j\tau_M\big),\qquad j=1,\ldots,M.
\end{equation*}
Similarly to Theorem~\ref{thm:existenceQVI}, we can verify that, for each fixed $\bfy_0\in\mathbb{R}^N$, the discretized quasi-variational inequality \eqref{eq:QVIdd} admits a solution $\mathbf{y}=\{\mathbf{y}_j^M\}_{j=1}^M$.  The {\em discrete version} of the optimal control problem $(\mathbb{P}_N)$ is formulated as follows: \vspace*{0.1in}

\textbf{Problem} $(\mathbb{P}_N^M)$. Given $\sigma>0$, a nonnegative mapping $\mathbf{f}:(0,T)\to\mathbb{R}^N$, and vectors $\mathbf{a},\;\mathbf{y}_0^{\mathrm{ref}}\in\mathbb{R}^N$, consider the discrete-time optimal control problem:
\begin{align*}
&\mathrm{minimize}\qquad J^M(\mathbf{y},\mathbf{y}_0):=\sum_{j=1}^M\tau_M\big\la\mathbf{a}, \mathbf{y}_j^M-\mathbf{y}_0\big\ra+\frac{\sigma}{2}\big\la\mathbf{y}_0-\mathbf{y}_0^{\mathrm{ref}},\mathbf{y}_0-\mathbf{y}_0^{\mathrm{ref}}\big\ra\\
&\text{over}\qquad\qquad \mathbf{y}_0^M,\mathbf{y}_1^M,\ldots,\mathbf{y}_M^M\in \mathbb{R}^N;\\
&\mathrm{subject \:\:to\:\:}\quad \mathbf{y}=\{\mathbf{y}_j^M\}_{j=1}^M\;\text{ solves }\; \mathrm{QVI}_N^M(\mathbf{y}_0),\\
&\hphantom{\mathrm{subject \:\: to \:\:}\quad}  \mathbf{y}_0\in\mathcal{A}.
\end{align*}
In this problem, the {\em dynamics constraints} can be written in the quasi-variational sweeping form
\begin{equation}\label{e:diff-incl}
\dot{\bfy}(t_j^M)\in F_j^M(\bfy(t_j^M),\bfy_0)\;\mbox{ for all }\;t^M_j\in(0,T)
\end{equation}
with $F_j^M(\mathbf{y},\bfy_0)$ from \eqref{mapF1}, the {\em control constraint} $\mathbf{y}_0\in\mathcal{A}$ is expressed in terms of the set
\begin{equation}\label{A}
\mathcal{A}:=\big\{\mathbf{z}\in\mathbb{R}^N\;\big|\;\mathbf{y}_0^{\mathrm{ref}}+\bm{\lambda}_0\le \mathbf{z}\le\mathbf{y}_0^{\mathrm{ref}}+\bm{\lambda}_1 \big\},
\end{equation}
where $\bm{\lambda}_0,\bm{\lambda}_1\in\mathbb{R}^N$ with $0\le\bm{\lambda}_0\le\blm_1$, and the {\em hidden state constraints} are given by 
\begin{equation}\label{e:ic}
-\big(M_\infty(\bfy(t_j^M),\bfy_0)\big)_i\le \big(\mathbf{D}_k\bfy(t_j^M)\big)_i\le \big(M_\infty(\bfy(t_j^M),\bfy_0)\big)_i
\end{equation}
with $i=1,\ldots,N$, $k=1,2$, and $j=1,\ldots,M$, where
the mapping $M_p$ is defined in \eqref{M_p}.\vspace*{0.08in}

Similarly to the proof of Theorem~\ref{thm:ExistPN}, we arrive at following existence theorem of optimal solutions. \vspace*{-0.1in}

\begin{thm}{\bf(existence of optimal solutions to discretized sweepings QVIs).}\label{thm:existencePnh} For each natural numbers $N$ and $M$, the discretized sweeping control problem $(\mathbb{P}_N^M)$ admits an optimal solution.
\end{thm}

It has been well understood in the developments of the discrete approximation method for Lipschitzian differential inclusions \cite{m95,m-book} and for sweeping control problems associated with variational inequalities \cite{ccmn,cm3,chhm1,chhm2} that 
optimal solutions to the discrete-time problems of the above type {\em strongly converge} in the suitable space topologies to the prescribed local minimizer of the original continuous-time problem. A similar result holds for the controlled quasi-variational sweeping process $(\mathbb{P}_N)$ and its discrete approximations $(\mathbb{P}_N^M)$ under consideration by imposing appropriate assumptions, while we postpone the precise clarification of this issue to our future research.\vspace*{0.05in}

Our further goal in this paper is to derive {\em necessary optimality conditions} for local minimizers of the discrete-time quasi-variational sweeping control problem $(\mathbb{P}^M_N)$ for each $N,M\in\N$. According to the previous discussions, such necessary optimality conditions for $(\mathbb{P}^M_N)$ can be viewed as {\em suboptimality} (almost optimality) condition for $(\mathbb{P}_N)$ and the original quasi-variational control problem $(\mathbb{P})$. 

Looking at the structure of each problem $({\mathbb P}^M_N)$ tells us that it can be reduced to a problem of {\em finite-dimensional optimization} while with a special type of (increasingly many) {\em geometric constraints} given in the unavoidably {\em nonconvex graphical} form induced by the very nature of the quasi-variational sweeping process. Handling such constraints require the usage of adequate tools of nonconvex variational analysis and generalized differentiation, which we briefly review in the next section.

\section{Generalized Differentiation for QVI Sweeping Dynamics}\label{sec:coderivative}

First we present here the generalized differential notions for sets, set-valued mappings, and extended-real-valued functions that are used in what follows. More details and references can be found in the books \cite{m-book,m18,rw}.

Following the geometric approach of \cite{m-book,m18}, we start with generalized normals to sets. Given a set $\Theta\subset\R^s$ locally closed around $\oz\in\Theta$, the (Mordukhovich, limiting) {\em normal cone} to $\Theta$ at $\oz$ is defined by
\begin{equation}\label{lim-nor}
N_\Theta(\oz):=\big\{v\in\R^s\;\big|\;\exists\,z_k\to\oz,\;w_k\in\Pi_\Theta(z_k),\;\al_k\ge 0\,\;\mbox{ with }\;\al_k(z_k-w_k)\to v\big\}, 
\end{equation}
where $\Pi_\Theta(z)$ stands for the (nonempty) Euclidean projector of $z\in\R^s$ onto $\Theta$. If $\Theta$ is convex, the normal cone \eqref{lim-nor} agrees with normal cone of convex analysis \eqref{nor-cone}, but otherwise \eqref{lim-nor} is {\em nonconvex} in very common situations, e.g., for the graph of $\ph(x):=|x|$ and the epigraph of $\ph(x):-|x|$ at $(0,0)\in\R^2$. Nevertheless, the normal cone \eqref{lim-nor} and the associated generalized differential constructions for mappings and functions defined below enjoy comprehensive {\em calculus rules} the proofs of which are based on the {\em variational/extremal principles} of variational analysis. 

Let $F\colon\R^n\tto\R^m$ be a set-valued mapping/multifunction with graph
$$
\gph F:=\big\{(x,y)\in\R^n\times\R^m\;\big|\;y\in F(x)\big\}
$$
locally closed around $(\ox,\oy)\in\gph F$. The {\em coderivative} of $F$ at $(\ox,\oy)$ {is defined via} the normal cone 
\eqref{lim-nor} to the graph of $F$ at this point by
\begin{equation}\label{e:cor}
D^*F(\ox,\oy)(u):=\big\{v\in\R^n\;\big|\;(v,-u)\in N_{{\rm\small gph}\,F}(\ox,\oy)\big\}\;\mbox{ for all }\;u\in\R^m.
\end{equation}
This is an extension to the case of nonsmooth and set-valued mappings the notion of the {\em adjoint operator} (matrix transposition) for the Jacobians $\nabla F(\ox)$ of single-valued smooth mappings in which case we have
$$
D^*F(\ox)(u)=\big\{\nabla F(\ox)^*u\big\},\quad u\in\R^m,
$$
where the indication of $\oy=F(\ox)$ is dropped in the coderivative notation.

Let $\ph\colon\R^n\to\oR:=(-\infty,\infty]$ be an extended-{real-valued} function that is lower semicontinuous (l.s.c.) around $\ox$ with $\ph(\ox)<\infty$, i.e., with $\ox\in\dom\ph$. Proceeding geometrically, the (first-order) {\em subdifferential} of the function $\ph$ at the point $\ox$ is defined as
\begin{equation}\label{e:sub}
\partial\ph(\ox):=\big\{v\in\R^n\;\big|\;(v,-1)\in N_{{\rm\small epi}\,\ph}\big(\ox,\ph(\ox)\big)\big\}
\end{equation}
via the normal cone to the epigraph $\epi\ph$ of $\ph$ at $(\ox,\ph(\ox))$ while observing that the subgradient mapping $\partial\ph$ admits various equivalent analytic descriptions that can be found in the aforementioned books. 

Following the ``dual derivative-of-derivative" scheme of \cite{m92}, we finally introduce the major second-order generalized differential construction used in the paper. Given $(\ox,\ov)\in\gph\partial\ph$ for an l.s.c.\ function $\ph\colon\R^n\to\oR$,
the {\em second-order subdifferential}, or the {\em generalized Hessian}, of $\ph$ at $\ox$ relative to $\ov$ is 
\begin{eqnarray}\label{2nd}
\partial^2\ph(\ox,\ov)(u):=\big(D^*\partial\ph\big)(\ox,\ov)(u),\quad u\in\R^n.
\end{eqnarray}
When $\ph$ is ${\cal C}^2$-smooth around $\ox$, we have the representation
\begin{equation*}
\partial^2\ph(\ox)(u)=\big\{\nabla^2\ph(\ox)u\big\}\;\mbox{ for all }\;u\in\R^n
\end{equation*}
via the (symmetric) Hessian matrix of $\ph$ at $\ox$. The well-developed {\em second-order calculus} is available for \eqref{2nd} in general settings, and explicit evaluations of this construction is given for major classes of functions important in applications to nonsmooth optimization, optimal control, and related topics; see, e.g., \cite{m-book,m18,mr} and the references therein. Note that coderivatives and second-order subdifferentials has been already used in \cite{mo} in the study of nondynamic finite-dimensional quasi-variational inequalities in the framework of generalized equations, which is totally different from our current consideration.\vspace*{0.03in}

To efficiently proceed in the setting of this paper, we modify $({\mathbb P}^M_N)$ a bit with replacing the constraint mapping $M_\infty$ in \eqref{e:ic} by its {\em smooth version} $\Tilde M_\infty$ for $p=\infty$. The corresponding set \eqref{setK}, with the replacement of $M_\infty$ by $\Tilde M_\infty$, is labeled as $\Tilde{\cal K}^\infty$. Define further
\begin{equation}\label{Theta}
\Theta:=\big\{(\bfy,\bfy_0)\in\R^N\times\R^N\big|\;g^l_{k}(\bfy,\bfy_0)\ge 0\big\},\quad\;l=1,\ldots,N,N+1,\ldots,2N,\quad\;k=1,2,
\end{equation}
via the twice continuously differentiable mapping $g\colon\R^{2N}\to\R^{4N}$ with the components
\begin{equation}\label{g}
g^i_{k}(\bfy,\bfy_0)=(\bfD_k\bfy)_i+\big(\tilde M_{\infty}(\bfy,\bfy_0)\big)_i,\quad g^{N+i}_k(\bfy,\bfy_0)=\big(\tilde M_{\infty}(\bfy,\bfy_0)\big)_i-(\bfD_k\bfy)_i,
\end{equation}
where $\bfy_i$ stands for the $i^{th}$ coordinate {of the underlying vector.}
\vspace{0.03in}

For our application to deriving necessary optimality conditions for problem $({\mathbb P}^M_N)$ with the {\em smoothed constraints} as above (no relabeling), we are going to compute the second-order subdifferential \eqref{2nd} of the indicator function $\ph:=\delta_\Theta(z)$ of the set $\Theta$ from \eqref{Theta}, i.e., such that $\dd_\Theta(z):=0$ if $z\in\Theta$ and $\dd_\Theta(z):=\infty$ otherwise. In this case, we have $\partial\ph=N_\Theta$ and $\partial^2\ph=D^*N_\Theta$. Recall that the {\em domain} (dom) of a set-valued mapping contains those points where the mapping has nonempty values,\vspace*{-0.05in}
  
\begin{thm}{\bf(second-order computation for the discretized QVI sweeping process).}\label{Th:co-cal} Consider problem $({\mathbb P}^M_N)$ with the smoothed constraints for any fix $N,M\in\N$, and let $F:=F_j^M$ be taken from \eqref{mapF1} with $p=\infty$ and with ${\cal K}^\infty$ replaced by $\tilde{\cal K}^\infty$, where $\bff_j^M$ is generated by $\bff$ in \eqref{f-discr}. Given $(\bfy,\bfy_0)\in\Theta$, assume that the gradient vectors $\{\nabla g^1_1(\bfy,\bfy_0),\ldots,\nabla g^{2N}_2(\bfy,\bfy_0)\}$ for the functions from \eqref{g} are linearly independent. Then there exists the collection of nonnegative multipliers $\lm_1,\ldots,\lm_{2N}$ uniquely determined by the equation
$-\nabla g(\bfy,\bfy_0)^*\lm=w+\bff$ for $\lm=(\lm_1,\ldots,\lm_{2N})$ such that 
\begin{eqnarray*}
\begin{array}{ll}
D^*F(\bfy,\bfy_0,w)(y)=\disp\bigcup_{\lm\ge 0,-\nabla g(\bfy,\bfy_0)\lm=w+\bff}\bigg\{\bigg(-\sum_{k=1}^2\sum^{2N}_{l=1}\lm^l_k\big\la\nabla^2_{\bfy}g^l_k(\bfy,\bfy_0),y\big\ra\}-\nabla_{\bfy}g(\bfy,\bfy_0)^*\gg,0\bigg)\bigg\}
\end{array}
\end{eqnarray*}
$\mbox{for all }\;y\in\dom D^*N_{\Tilde\K^{\infty}(\bfy,\bfy_0)}\big(\bfy,w+\bff\big)$, where the coderivative domain is given by
\begin{equation*}
\begin{aligned}
\dom D^*N_{\tilde\K^{\infty}(\bfy,\bfy_0)}(\bfy,w+\bff)=&\big\{y\,\big|\;\exists\,\lm\ge0\;\mbox{ such that }\;-\nabla g(\bfy,\bfy_0)\lm=w+\bff,\\
&\lm^l_k\la\nabla g^l_k(\bfy,\bfy_0),y\ra=0\;\mbox{ for }\;l=1,\ldots,2N,\;k=1,2\big\}
\end{aligned}
\end{equation*}
with $\gg^l_k=0$ if either $g^l_k(\bfy,\bfy_0)>0$ or $\lm^l_k=0$ and $\la\nabla g^l_k(\bfy,\bfy_0),y\ra>0$, and with $\gg^l_k\ge 0$ if $g^l_k(\bfy,\bfy_0)=0,\;\lm^l_k=0$, and $\la\nabla g^l_k(\bfy,\bfy_0),y\ra<0$.
\end{thm}
{\bf Proof}. Define the set-valued mapping $G$ and the single-valued smooth mapping $\Tilde\bff$ by, respectively,
\begin{equation*}
G(\bfy,\bfy_0):=N_{\tilde\K^{\infty}(\bfy,\bfy_0)}(\bfy)\;\mbox{ and }\;\Tilde\bff(\bfy,\bfy_0):=\bff. 
\end{equation*}
The coderivative sum rule from \cite[Theorem~1.62]{m-book} tells us that
\begin{equation*}
z^*\in\nabla\Tilde\bff(\bfy,\bfy_0)^*y+D^*G\big(\bfy,\bfy_0,w+\bff\big)(y)
\end{equation*}
for any $y\in\dom D^*N_{\Tilde\K^{\infty}(\bfy,\bfy_0)}(\bfy,w+\bff)$ and $z^*\in D^*F(\bfy,\bfy_0,w)(y)$. Observe further that
\begin{equation*}
G(\bfy,\bfy_0)=N_{\Tilde\K^{\infty}(\bfy,\bfy_0)}\circ\Tilde g(\bfy,\bfy_0)\;\mbox{ with }\;\Tilde g(\bfy,\bfy_0):=\bff,
\end{equation*}
where the Jacobian of latter mapping is obviously of full rank. Employing further the coderivative chain rule from \cite[Theorem~1.66]{m-book} to the above composition for $G$ yields
\begin{equation}\label{e:cf}
z^*\in\nabla\Tilde\bff(\bfy,\bfy_0)^*y+\nabla\Tilde g(\bfy,\bfy_0)^*D^*N_{\Tilde\K^{\infty}(\bfy,\bfy_0)}\big(\bfy,w+\bff\big)(y).
\end{equation}
To deduce finally from \eqref{e:cf} the exact formulas claimed in the theorem, we use for representing $D^*N_{\Tilde\K^{\infty}}$ the second-order calculation for inequality constraint systems
taken from \cite[Theorem~3.3]{hos} in the case of the linear independence condition imposed in this theorem.$\h$

\section{Necessary Optimality Conditions for Discrete-Time Problems}\label{sec:NecOptCond} The main result of this section provides constructive necessary optimality conditions for each problem $({\mathbb P}^M_N)$ expressed in terms of its initial data. As discussed above, for $N,M\in\N$ sufficiently large such conditions may be viewed as {\em suboptimality} conditions for problems with the semi-discrete $({\mathbb P}_N)$ and continuous-time $({\mathbb P})$ dynamics.

To accomplish our goal with taking into account the complexity of smoothed problem $({\mathbb P}^M_N)$, we split the derivation of necessary optimality conditions into two theorems. The first theorem presents necessary optimality conditions for $({\mathbb P}^M_N)$ that involve the limiting normal cone to graphs of discrete velocity mappings, i.e., their {\em coderivatives}. The final result benefits from the coderivative computations for such mappings furnished in Theorem~\ref{Th:co-cal} and thus provides necessary optimality conditions for $({\mathbb P}^M_N)$ explicitly expressed in terms of the problem data.\vspace*{0.05in}

Our general scheme of deriving necessary optimality conditions for $({\mathbb P}^M_N)$ is similar to the one in \cite{cm3} addressed an optimal control problem for a sweeping process over state-independent and canonically controlled {\em prox-regular} moving sets of the type
\begin{equation}\label{prox}
C(t)=C+u(t)\;\mbox{ with }\;C:=\big\{x\in\R^n\;\big|\;g_i(x)\ge 0,\;i=1,\ldots,m\big\},
\end{equation}
where $x$ and $u$ stand for the state and control variables, respectively. However, the setting of problem $({\mathbb P}^M_N)$ is very different from \cite{cm3}. First of all, we have the {\em state-dependent} moving sets (the essence of QVI) with nonlinear control functions. Indeed, the counterpart of $C(t)$ in \eqref{prox} is the set $\Tilde{\cal K}^\infty(\bfy,\bfy_0)$ depending on both state $\bfy$ and control $\bfy_0$ variables being described in form \eqref{Theta} via the functions $g^i_k(\bfy,\bfy_0)$ from \eqref{g}. Observe also that, in contrast to \eqref{prox}, the functions from \eqref{g} are ${\cal C}^2$-smooth while may be {\em nonconvex}, which does not allow us the claim the prox-regularity of the moving sets in $({\mathbb P}^M_N)$ as in \cite{cm3}. Nevertheless, we can proceed with deriving necessary optimality conditions for problem $({\mathbb P}^M_N)$ by reducing it to a problem of {\em mathematical programming} with functional and geometric constraints and then using the machinery of {\em variational analysis and generalized differentiation} discussed above.\vspace*{0.05in}

Here is the first theorem involving coderivatives (without their explicit computations) of the mappings in the {\em smoothed} dynamic constraints \eqref{e:diff-incl} with $F^M_j$  defined by
\begin{equation}\label{mapF2}
F_j^M(\mathbf{y}^M_j,\bfy_0):=-N_{\tilde\K^{\infty}(\mathbf{y}^M_j,\bfy_0)}(\mathbf{y}^M_j)+\bff_j^M,\quad j=1,\ldots,M,
\end{equation}
according to our previous discussions, where the state-control dependent moving sets $\tilde\K^{\infty}(\mathbf{y}^M_j,\bfy_0)$ are generated by the functions $g^i_k$ from \eqref{g}.\vspace*{-0.05in}

\begin{thm}{\bf(coderivative-based necessary optimality conditions for discretized QVI problems).}\label{NOC}
Let $\(\obfy^M,\obfy_0\)=(\obfy_1^M,\ldots,\obfy_{M}^M,\obfy_0)$ be an optimal solution to problem $(\P^M_N)$ with smoothed constraints, and let $F_j:=F^M_j$ be taken from \eqref{mapF2}. Assume that the  gradients $\{\nabla g^1_1(\obfy,\obfy_0),\ldots,\nabla g^{2N}_2(\obfy,\obfy_0)\}$ are linearly independent. Then there exist dual elements $\lm^M\ge 0$, $\al^{kM}=\(\al_{1}^{kM},\ldots,\al_{2N}^{kM}\)\in\R^{2N}_+$, and $p^M_j\in\R^N$ as $j=1,\ldots,M$ satisfying the conditions
\begin{equation}\label{NOC1}
\lm^M+\|\al^{kM}\|+\sum_{j=1}^{M}\n p_{j}^{M}\en+\n\psi\en\ne 0,
\end{equation}
\begin{equation}\label{con:al1}
\al_{l}^{kM}g^l_{k}(\obfy^M_M,\obfy_0)=0\;\mbox{ for all }\;\;l=1,\ldots,2N\;\mbox{ and }\;k=1,2,
\end{equation}
\begin{equation}\label{pNN}
p^{M}_{M}=\sum_{k=1}^2\sum^{2N}_{l=1}\al_{l}^{kM}\nabla_{\obfy^M_M} g^l_k(\obfy^M_M,\obfy_0),
\end{equation}
\begin{equation}\label{inclu}
\begin{aligned}
\Bigg(\frac{p_{j+1}^{M}-p_{j}^{M}}{\tau_M}-\lm^M\mathbf{a}^{\intercal} ,-\frac{1}{\tau_M}\lm^M\(-\tau_m\mathbf{a}^{\intercal} +\sigma\obfy_0\)+\dfrac{1}{\tau_M}\sum_{k=1}^2\sum^{2N}_{l=1}\al_{l}^{kM}\nabla_{\obfy_0}g^l_k(\obfy^M_M,\obfy_0),p_{j+1}^{M}\Bigg)\\
\in\(0,\disp\frac{1}{\tau_M}\psi,0\)+N\(\(\obfy_j^M,\obfy_0,-\dfrac{\obfy_{j+1}^M-
\obfy_{j}^M}{\tau_M}\);\gph F_j\)
\end{aligned}
\end{equation}
for all $\,j=1,\ldots,M-1$ and $k=1,2$ together with the inclusion
\begin{equation}\label{psi}
\psi\in N_{\mathcal{A}}\(\obfy_0\),
\end{equation}
where the functions $g^i_k$ and the set ${\cal A}$ are taken from \eqref{g} and \eqref{A}, respectively.
\end{thm}
{\bf Proof}. Fix $\ve>0$ and consider the vector
\begin{equation*}
z:=\(\bfy_1^M,\ldots,\bfy_{M}^M,\bfy_0,\bfY^M_1,\ldots,\bfY^M_{M-1}\).
\end{equation*}
Then $(\P^M_N)$ is equivalent to problem of mathematical programming $(MP)$ with respect $z$:
\begin{equation*}
\begin{aligned}
\textrm{minimize }\phi_0(z):=\sum_{j=0}^{M-1}\int_{t_{j}^M}^{t_{j+1}^M}\big\la\mathbf{a}, \bfy_j^M-\bfy_0\big\ra\mathrm{d}t+\frac{\sigma}{2}\big\la\bfy_0-\bfy_0^{\mathrm{ref}},\bfy_0-\mathbf{y}_0^{\mathrm{ref}}\big\ra
\end{aligned}
\end{equation*}
subject to the functional and geometric constraints
\begin{equation*}
H_j(z):=\bfy_{j+1}^M-\bfy_{j}^M-\tau_M\bfY_j^M=0,\;j=0,\ldots,M-1,
\end{equation*}
\begin{equation*}
L^l_k(z:)=-g^{l}_k(\bfy^M_M,\bfy_0)\le 0,\;\;l=1,\ldots,2N,\;\;k=1,2,
\end{equation*}
\begin{equation*}
z\in\Xi_j:=\left\{z\;\bigg{|}\;\bfY_{j}^M\in F_{j}\(\bfy_{j}^M,\bfy_0\)\right\},\;j=1,\ldots,M,
\end{equation*}
\begin{equation*}
z\in\O:=\left\{z\;\bigg{|}\;\bfy_0\in\mathcal{A}\right\}.
\end{equation*}
Applying now the necessary optimality conditions from \cite[Theorem~6.5]{m18} to finite-dimensional mathematical programming problem $(MP)$ at its optimal solution
\begin{eqnarray*}
\oz:=\(\obfy_1^M,\ldots,\obfy^{M}_M,\obfy_0,\obfY_1^M,\ldots,\obfY_{M-1}^M\)\in\R^{2MN}
\end{eqnarray*}
gives us dual elements $\lm^M\ge 0$, $p^M_j\in\R^N,\;j=2,\ldots,M$, $\al^{kM}=(\al_{1}^{kM},\ldots,\al_{2N}^{kM})\in\R^{2N}_+$ as $k=1,2$, and
\begin{eqnarray*}
z^*_j=\(\bfy^*_{1j},\ldots,\bfy^*_{Mj},\bfy_0^*,\bfY^*_{1j},\ldots,\bfY^*_{(M-1)\,j}\)\in\R^{2MN}
\end{eqnarray*}
for $j=1,\ldots,M$, which are not zero simultaneously, while satisfying the following relationships:
\begin{equation}\label{sncs}
z^*_j\in N_{\Xi_j}(\oz)+N_{\O}(\oz),\quad j\in\{1,\ldots,M,
\end{equation}
\begin{equation}\label{main}
-z^*_1-\ldots-z^*_{M}\in\lm^M\partial\phi_0(\oz)+\sum_{k=1}^2\sum^{2N}_{l=1}\al_{l}^{kM}\nabla L^l_{k}(\oz)+\sum_{j=0}^{M-1}\nabla H_j(\oz)^*p_{j+1}^{M},
\end{equation}
\begin{equation}\label{eq:al1}
\al_{l}^{kM}L^l_k\(\oz\)=0\;\mbox{ as }\;l=1,\ldots,2N\;\mbox{ and }\;k=1,2.
\end{equation}
To specify more, note that in \eqref{sncs} we apply the normal cone intersection formula from \cite[Theorem~2.16]{m18} to $\oz\in\O_j\cap\Xi_j$ for $j=1,\ldots,M-1$, where the qualification condition therein holds due to the graphical structure of the sets $\Xi_j$ and the coderivative computation from Theorem~\ref{Th:co-cal}. Furthermore, the structure of the sets $\O_j$ and $\Xi_j$ together with \eqref{sncs} leads us to the relationships
$$\(\bfy^\ast_{11},\ldots,\bfy^\ast_{M1},\bfy^\ast_{0},\bfY^\ast_{11},\ldots,\bfY^\ast_{(M-1)1}\)\in N_{\Xi_1}\(\oz\)+N_{\Omega}(\oz),$$
$$\(\bfy^\ast_{12},\ldots,\bfy^\ast_{M2},\bfy^\ast_{0},\bfY^\ast_{12},\ldots,\bfY^\ast_{(M-1)2}\)\in N_{\Xi_2}\(\oz\)+N_{\Omega}(\oz),$$
$$\ldots$$
$$\(\bfy^\ast_{1M},\ldots,\bfy^\ast_{MM},\bfy^\ast_{0},\bfY^\ast_{1M},\ldots,\bfY^\ast_{(M-1)M}\)\in N_{\Xi_M}\(\oz\)+N_{\Omega}(\oz).$$
In this way we arrive at the inequality
\begin{equation*}
\begin{array}{ll}
\bigg\la\(\bfy^\ast_{11},\ldots,\bfy^\ast_{M1},\bfy^\ast_0,\bfY^\ast_{11},\ldots,\bfY^\ast_{(M-1)1}\),\Big(\big(\obfy_1,\ldots,\obfy_M,\obfy_0,\obfY_1\ldots,\obfY_{M-1}\Big)\\
\qquad\qquad\qquad\qquad\qquad\qquad\qquad\qquad\qquad\quad-\Big(\bfy_1,\ldots,\bfy_M,\bfy_0,\bfY_1,\ldots,\bfY_{M-1}\Big)\bigg)\bigg\ra\le 0,
\end{array}
\end{equation*}
where $\obfY_j\in F_j(\obfy_j,\obfy_0)$ and $\bfY_j\in F_j(\bfy_j,\bfy_0)$ for all $j=1,\ldots,M-1$, and where $\obfy_0,\bfy_0\in\mathcal{A}$. Combining the above verifies that $\bfy^\ast_{ij}=\bfY^\ast_{ij}=0$ if $i\ne j$, for all $j=1,\ldots,M$. The obtained relationships ensure that the inclusions in \eqref{sncs} are equivalent to
\begin{equation}\label{e:5.18*}
\(\bfy^*_{jj},\bfy_0^*-\psi,-\bfY^*_{jj}\)\in N\(\(\obfy_j^M,\obfy_0,-\dfrac{\obfy_{j+1}^M-\obfy_j^M}{\tau_M}\);\gph F_j\),\quad j=1,\ldots,M-1,
\end{equation}
while all the other components of $z^*_j$ are equal to zero for $j=1,\ldots,M-1$. We also get from above that $\psi\in N\(\obfy_0;\mathcal{A}\)$, which justifies \eqref{psi}. It follows furthermore that
\begin{eqnarray*}
-z^*_1-\ldots-z^*_{M}=\big(-\bfy^*_{11},\ldots,-\bfy^*_{M-1\,M-1},0,-\bfy_0^*,-\bfY^*_{11},\ldots,-\bfY^*_{M-1\,M-1}\big).
\end{eqnarray*}
The set on the right-hand side of \eqref{main} is represented by
\begin{eqnarray*}
\lm^M\partial\phi_0(\oz)+\sum_{k=1}^2\sum^{2N}_{l=1}\al_{l}^{kM}\nabla L^{l}_{k}(\oz)+\sum_{j=0}^{M-1}\nabla H_j(\oz)^*p_{j+1}^{M}.
\end{eqnarray*}
Using the definitions of $L^l_k$ and $H_j$, we easily obtain the equality
$$
\(\sum_{k=1}^2\sum^{2N}_{l=1}\al_{l}^{kM}\nabla L^l_{k}(\oz)\)_{(\obfy_j^M,\obfy_0,\obfY_j^M)}=\(-\sum_{k=1}^2\sum^{2N}_{l=1}\al_{l}^{kM}\nabla_{\obfy_j^M} g^l_k(\obfy^M_M,\obfy_0),-\sum_{k=1}^2\sum^{2N}_{l=1}\al_{l}^{kM}\nabla_{\obfy_0} g^l_k(\obfy^M_M,\obfy_0),0\)
$$
for $j=1,\ldots,M-1$ together with the representations
$$
\(\sum_{j=0}^{M-1}\nabla H_j(\oz)^*p_{j+1}^M\)_{\obfy_j^M}=
\(-p_{1}^{M},p_{1}^{M}-p_{2}^{M},\ldots,p_{j}^{M}-p_{j+1}^{M},\ldots,p_{M-1}^{M}-p_{M}^{M},p^{M}_{M}\),
$$
$$
\(\sum_{j=0}^{M-1}\nabla H_j(\oz)^*p_{j+1}^M\)_{\obfY_j^M}=\big(-\tau_Mp_{1}^{M},-\tau_Mp_{2}^{M},\ldots,-\tau_Mp^{M}_{M}\big).
$$
The set $\lm^M\partial\phi_0(\oz)$ is represented as the collection of
$$
\lm^M\((\tau_m\mathbf{a} )_{\bfy_1^M},\ldots,(\tau_m\mathbf{a} )_{\bfy^M_{M-1}},0,(-\tau_m\mathbf{a} +\sigma \obfy_0)_{(\bfy_0)},0_{\bfY_1^M},\ldots,0_{\bfY_{M-1}^M}\).
$$
Combining the above gives us the relationships
\begin{equation}\label{ey1}
-\bfy^*_{11}=\lm^M\tau_\mathbf{a}-p_{2}^{M},
\end{equation}
\begin{equation}\label{ey2}
-\bfy^*_{jj}=\lm^M\tau_m\mathbf{a}+p^{j}_{M}-p^{j+1}_{M},\;\;j=2,\ldots,M-1,
\end{equation}
\begin{equation}\label{ey3}
0=p_{M}^{M}-\sum_{k=1}^2\sum^{2N}_{l=1}\al_{l}^{kM}\nabla_{\obfy^M_M} g^l_k(\obfy^M_M,\obfy_0),
\end{equation}
\begin{equation}\label{ey0}
-\bfy_0^*=\lm^M\(-\tau_m\mathbf{a} +\sigma\obfy_0\)-\sum_{k=1}^2\sum^{2N}_{l=1}\al_{l}^{kM}\nabla_{\obfy_0} g^l_k(\obfy^M_M,\obfy_0),
\end{equation}
\begin{equation}\label{eY}
-\bfY^*_{jj}=-\tau_Mp_{j+1}^{M},\;\;j=0,\ldots,
M-1.
\end{equation}
Using the obtained representations, we can now proceed with completing the proof of the theorem. First observe that the transversality condition \eqref{pNN} follows directly from \eqref{ey3}. Next we extend the vector $p^M$ by adding the component $p_1^M:=\bfy^*_{1M}$. This tells us by \eqref{ey2}, \eqref{ey0}, and \eqref{eY} that
$$
\begin{aligned}
\frac{\bfy^*_{jj}}{\tau_M}&=\dfrac{p_{j+1}^{M}-p_{j}^{M}}{\tau_M}-\frac{l\lm^M \mathbf{a}}{\tau_M},\\
\frac{\bfy_0^*}{\tau_M}&=-\dfrac{1}{\tau_M}\lm^M\(-\tau_m\mathbf{a}+\sigma \obfy_0\)+\dfrac{1}{\tau_M}\sum_{k=1}^2\sum^{2N}_{l=1}\al_{l}^{kM}\nabla_{\obfy_0} g^l_k(\obfy^M_M,\obfy_0),\\
\frac{\bfY^*_{jj}}{\tau_M}&=p_{j+1}^{M}.
\end{aligned}
$$
Substituting these relationships into the left-hand side of \eqref{e:5.18*} and taking into account the equalities obtained in  \eqref{eq:al1}, \eqref{ey2}, \eqref{ey0}, and \eqref{eY} verify the optimality conditions claimed in  \eqref{con:al1}--\eqref{inclu}.

It remains to justify the nontriviality condition \eqref{NOC1}. On the contrary, suppose that $\lm^M=0,\,\psi=0,\,\al^{kM}=0$, and $p_{j}^{M}=0$ for all $j=1,\ldots,M-1$, which implies in turn that $\bfy^*_{1M}=p_{1}^{M}=0$. Then we deduce from \eqref{ey3} that $p^{M}_{M}=0$, and so $p_{j}^{M}=0$ for all $j=1,\ldots,M$. It follows from \eqref{ey1}, \eqref{ey2}, and \eqref{ey0} that $\(\bfy^*_{jj},\bfy_0^*\)=0$ for all $j=1,\ldots,M-1$. By \eqref{eY} we have that $\bfY^*_{jj}=0$ whenever $j=1,\ldots,M-1$. Since all the components of $z^*_j$ different from $(\bfy^*_{jj},\bfy_0^*,\bfY^*_{jj})$ are obviously zero for $j=1,\ldots,M-1$, this tells us that $z^*_{j}=0$ for such $j$. Employing $\bfy^*_{1M}=p_{1}^{M}=0$ ensures that $z^*_{M}=0$ while the other components of this vector are zero. Overall, $z^*_j=0$ for all $j=1,\ldots,M$, and thus the nontriviality condition for $(MP)$ fails. The obtained contradiction completes the proof. $\h$\vspace*{0.05in}

The final result of this paper establishes necessary optimality conditions for smoothed problem $({\mathbb P}^M_N)$ expressed entirely in terms of the initial problem data. The desired conditions are derived by incorporating the second-order calculations of Theorem~\ref{Th:co-cal} into the corresponding conditions of Theorem~\ref{NOC} in the case where the mappings $F_j=F^M_j$ therein are given by \eqref{mapF2}.\vspace*{-0.05in}

\begin{thm}{\bf(explicit necessary conditions for discretized QVI sweeping control problems).}\label{Th:OC-DP} Let $\oz^M=(\obfy^M,\obfy_0)$ be an optimal control $\oz^M=(\obfy^M,\obfy_0)$ to the smoothed problem $(\P^M_N)$ with the sweeping dynamics defined by \eqref{mapF2} under the assumptions of Theorem~{\rm\ref{NOC}}. Then there exist dual elements $(\lm^M,\al^{kM},p^M)$ and $\psi\in N_{\mathcal{A}}\(\obfy_0\)$ together with vectors $\eta^{kM}_j=\(\eta^{kM}_{1j},\ldots,\eta^{kM}_{2Nj}\)\in\R^{2N}_+$ as $j=1,\ldots,M,\,k=1,2$ and $\gg^{kM}_j=\(\gg^{kM}_{1j},\ldots,\gg^{kM}_{2Nj}\)\in\R^{2N}$ as $j=1,\ldots,M-1$ and $k=1,2$ such that the following relationships hold:\\[1ex]
$\bullet$ {\sc nontriviality condition}
\begin{equation}\label{e:dac26}
\lm^M+\|\eta_M^{kM}\|+\sum^{M}_{j=1}\|p^{M}_j\|\not=0.
\end{equation}
$\bullet$ {\sc dynamic relationships} for all $j=1,\ldots,M-1$:
\begin{equation}\label{e:dac15}
\dfrac{\obfy^M_{j+1}-\obfy^M_j}{\tau_M}+\bff_j^M=-\sum_{k=1}^2\sum_{l\in I(\obfy^M_j)}\eta^{kM}_{ij}\nabla_{\obfy^M_j} g^l_k(\obfy^M_j,\obfy_0),
\end{equation}
\begin{equation}\label{e:dac16}
\dfrac{p^{M}_{j+1}-p^{M}_j}{\tau_M}-\frac{\lm^MT\mathbf{a}^{\intercal} }{\tau_M}=-\sum_{k=1}^2\sum^{2N}_{l=1}\eta^{kM}_{lj}\bigg\la\nabla^2_{\obfy^M_j}g^l_k(\obfy^M_j,\obfy_0),-p^{M}_{j+1}\bigg\ra-\sum_{k=1}^2\sum^{2N}_{l=1}\gg^{kM}_{lj}\nabla_{\obfy^M_j} g^l_k(\obfy^M_j,\obfy_0),
\end{equation}
\begin{equation}\label{e:dac17}
-\dfrac{1}{\tau_M}\lm^M\(T\mathbf{a}^{\intercal}+\sigma \obfy_0\)+\dfrac{1}{\tau_M}\sum_{k=1}^2\sum^{2N}_{l=1}\eta_{lM}^{kM}\nabla_{\obfy_0} g^l_k(\obfy^M_M,\obfy_0)-\dfrac{1}{\tau_M}\psi=0.
\end{equation}
$\bullet$ {\sc transversality condition}
\begin{equation}\label{e:dac19}
p^{M}_{M}=-\lm^MT\mathbf{a}+\sum_{k=1}^2\sum^{2N}_{l=1}\eta_{lM}^{kM}\nabla_{\obfy^M_M} g^l_k(\obfy^M_M,\obfy_0).
\end{equation}
$\bullet$ {\sc complementarity slackness conditions}
\begin{equation}
\label{e:dac20}
g^l_k(\obfy^M_j,\obfy_0)>0\Lto\eta^{kM}_{lj}=0,
\end{equation}
\begin{equation}
\label{e:dac21}
\big[g^l_k(\obfy^M_j,\obfy_0)>0\mbox{ or }\eta^{kM}_{lj}=0,\;\la \nabla g^l_k(\obfy^M_j,\obfy_0),-p^{M}_{j+1}\ra>0\big]\Lto\gg^{kM}_{lj}=0,
\end{equation}
\begin{equation}
\label{e:dac22}
\big[g^l_k(\obfy^M_j,\obfy_0)=0,\,\eta^{kM}_{lj}=0,\mbox{ and }
\la\nabla g^l_k(\obfy_j^M,\obfy_0),-p^{M}_{j+1}\ra<0\big]\Lto\gg^{kM}_{lj}\ge 0
\end{equation}
for $j=1,\ldots,M-1$, $l=1,\ldots,2N$, and $k=1,2$. Furthermore, we have the  implications
\begin{equation}
\label{e:dac23}
g^l_k(\obfy^M_j,\obfy_0)>0\Lto\gg^{kM}_{lj}=0\;\mbox{for}\;\;j=1,\ldots,M-1,\;\;l=1,\ldots,2N,\;\mbox{ and }\;k=1,2,
\end{equation}
\begin{equation}\label{e:dac24}
g^l_k(\obfy^M_M,\obfy_0)>0\Lto\eta^{kM}_{lM}=0\;\;\mbox{for}\;\;l=1,\ldots,2N\;\mbox{ and}\;\; k=1,2,
\end{equation}
\begin{equation}\label{e:dac25}
\eta^{kM}_{lj}>0\Lto\la\nabla g^l_k(\obfy^M_j,\obfy_0),-p^{M}_{j+1}\ra=0.
\end{equation}
\end{thm}
{\bf Proof}. The adjoint dynamic inclusion \eqref{inclu} of Theorem~\ref{NOC} can be written by the coderivative definition \eqref{e:cor} in the coderivative inclusion form
\begin{equation}\label{e:dac27}
\begin{aligned}
&\bigg(\frac{p_{j+1}^{M}-p_{j}^{M}}{\tau_M}-\frac{\lm^MT\mathbf{a} }{\tau_M},-\frac{1}{\tau_M}\lm^M\(T\mathbf{a}^{\intercal} +\sigma \obfy_0\)+\dfrac{1}{\tau_M}\sum_{k=1}^2\sum^{2N}_{l=1}\al_{l}^{kM}\nabla_{\obfy_0} g^l_k(\obfy^M_M,\obfy_0)-\frac{1}{\tau_M}\psi\bigg)\\
&\in D^*F_j\(\obfy_j^M,\obfy_0,-\dfrac{\obfy_{j+1}^M-\obfy_j^M}{\tau_M}\)(-p_{j+1}^{M}),\;\;j=1,\ldots,M-1.
\end{aligned}
\end{equation}
It follows from \eqref{mapF2} and the inclusions $\dfrac{\obfy^M_{j+1}-\obfy^M_j}{-\tau_M}-\bff_j\in N(\obfy^M_j;\Tilde\K^{\infty}(\obfy_j^M,\bfy_0))$ for $j=1,\ldots,M-1$ that the exist vectors $\eta^{kM}_j\in\R^{4N}_+$ as $j=1,\ldots,M-1$ and $k=1,2$ such that the conditions in \eqref{e:dac15} and \eqref{e:dac20} are satisfied. Employing the second-order formula from Theorem~\ref{Th:co-cal} with $\bfy:=\obfy^M_j,\,\bfy_0:=\obfy_0,\,w:=\dfrac{\obfy^M_{j+1}-\obfy^M_j}{-\tau_M}$, and $y:=-p^{M}_{j+1}$ and combining this with the domain formula therein give us vectors $\gg^{kM}_j\in\R^{4N}$ for which we have the equalities
$$
\begin{aligned}
&\bigg(\frac{p_{j+1}^{M}-p_{j}^{M}}{\tau_M}-\frac{\lm^MT\mathbf{a} }{\tau_M},-\frac{1}{\tau_M}\lm^M\sigma+\dfrac{1}{\tau_M}\sum_{k=1}^2\sum^{2N}_{l=1}\al_{l}^{kM}\nabla_{\obfy_0} g^l_k(\obfy^M_M,\obfy_0)-\frac{1}{\tau_M}\psi\bigg)\\
=&\bigg(-\sum_{k=1}^2\sum^{2N}_{l=1}\eta^{kM}_{lj}\bigg\la\nabla^2_{\obfy^M_j}g^l_k(\obfy^M_j,\obfy_0),-p^{M}_{j+1}\bigg\ra-\sum_{k=1}^2\sum^{2N}_{l=1}\gg^{kM}_{lj}\nabla_{\obfy^M_j} g^l_k(\obfy_j^M,\obfy_0),0\bigg)
\end{aligned}
$$
whenever $j=1,\ldots,M-1$. This clearly ensures the fulfillment of all the conditions claimed in \eqref{e:dac16}, \eqref{e:dac17}, \eqref{e:dac21}, and \eqref{e:dac22}. Now we denote $\eta_M^{kM}:=\al^{kM}$, where $\al^{kM}$ are taken from Theorem~\ref{NOC}, and note that $\eta^{kM}_j\in\R^{4N}_+$ for $j=1,\ldots,M$. Thus we get \eqref{NOC1} and deduce the transversality condition \eqref{e:dac19} from \eqref{pNN}. Observe also that \eqref{e:dac24} follows immediately from \eqref{con:al1} and the construction of $\eta_M^{kM}$, and that the adjoint inclusion \eqref{e:dac27} readily yields 
$$
-p^{M}_{j+1}\in\dom D^*N_{\Tilde\K^{\infty}(\obfy^M_j,\bfy_0)}\bigg(\obfy^M_j,\dfrac{\obfy^M_{j+1}-\obfy^M_j}{-\tau_M}+\bff_j\bigg).
$$
Based on this and coderivative formula from Theorem~\ref{Th:co-cal}, it is easy to check that  \eqref{e:dac25} is satisfied.

It remains to verify the nontriviality condition \eqref{e:dac26} taking into account the imposed gradient linear independence condition. On the contrary, suppose that \eqref{e:dac26} is violated, i.e., $\lm^M=0,\;\eta_{lM}^{kM}=0$ for $l=0,\ldots,2N,\,k=1,2$, and that $p^{M}_{j}=0$ for $j=1,\ldots,M$. Then it follows from \eqref{e:dac19} with 
$$
\disp\sum_{k=1}^2\sum^{2N}_{l=1}\eta^{kM}_{lM}\nabla g^l_k(\obfy^M_M,\obfy_0)=0
$$
that $p_{M}^M=0$. Then \eqref{e:dac17} tells us that $\psi\equiv 0$, and hence \eqref{e:dac16} implies that
$$
\disp\sum_{k=1}^2\sum^{2N}_{l=1}\gg^{kM}_{lj}\nabla g^l_k(\obfy^M_j,\obfy_0)=0\;\mbox{ for }\; j=1,\ldots,M-1.
$$
This contradicts the fulfillment of \eqref{NOC1} and thus verifies \eqref{e:dac26}. The proof is complete.$\h$

\section{Concluding Remarks}\label{conclusion}

This paper is the first attempt to study optimal control problems governed by evolutionary quasi-variational inequalities of the parabolic type that arise in the formation and growth modeling of granular cohensionless material. The formulated mathematical problem is revealed to be very complicated due to the presence of nonsmooth and nonconvex {gradient constraints} and thus calls for developing various regularization and approximation procedures for its efficient investigation and solution. Designing such procedures and verifying their well-posedness, we arrive at an adequate version described as optimal control of a discrete-time quasi-variational sweeping process, which is different from those previously considered in the literature.
Nevertheless, employing powerful tools of variational analysis and generalized differentiation brings us to the collection of necessary optimality conditions expressed entirely via the initial data of the original problem. These conditions are derived in Theorem~\ref{Th:OC-DP}.

{Some future research directions include, designing} efficient numerical algorithms {for the} system of optimality conditions presented in {Theorem~\ref{Th:OC-DP}. 
More work is needed to establish convergence analysis of some of the  regularization and approximation procedures constructed in this paper, for instance, regularization of the gradient constraints. This will, in particular, be critical to obtain optimality conditions for} $({\mathbb P}_N)$ and $({\mathbb P})$ by passing to the limit from those established in Theorem~\ref{Th:OC-DP} {for the fully discrete problem}.

\end{document}